\documentclass[12pt]{article}
\usepackage{amssymb}
\usepackage{amsmath}
\usepackage{amsfonts}

\newcommand{\be}{\begin{equation}}
\newcommand{\ee}{\end{equation}}
\newcommand{\bea}{\begin{eqnarray*}}
\newcommand{\eea}{\end{eqnarray*}}
\newcommand{\ba}{\begin{array}}
\newcommand{\ea}{\end{array}}
\newcommand{\bi}{\begin{itemize}}
\newcommand{\ei}{\end{itemize}}
\newcommand{\bc}{\begin{center}}
\newcommand{\ec}{\end{center}}
\newcommand{\bfr}{\begin{flushright}}
\newcommand{\efr}{\end{flushright}}

\input{tcilatex}

\begin{document}

\title{Semigroups of Operators on Spaces of Fuzzy-Number-Valued Functions
with Applications to Fuzzy Differential Equations}
\author{Ciprian G. Gal \\
Department of Mathematics\\
Florida international University\\
Miami, FL, 33199, E-mail: cgal@fiu.edu \\
\\
Sorin G. Gal \\
Department of Mathematics\\
University of Oradea, Romania\\
3700 Oradea, Romania\\
E-mail: galso@uoradea.ro}
\date{}
\maketitle

\begin{abstract}
In this paper we introduce and study semigroups of operators on spaces of
fuzzy-number-valued functions, and various applications to fuzzy
differential equations are presented. Starting from the space of fuzzy
numbers $\mathbb{R}_{F}$, many new spaces sharing the same properties are
introduced, as for example, with similar notations as in classical
functional analysis: $C([a,b];\mathbb{R}_{F})$, $C^{p}([a,b];\mathbb{R}_{F})$%
, $L^{p}([a,b];\mathbb{R}_{F})$, and so on. We derive basic operator theory
results on these spaces and new results in the theory of semigroups of
linear operators on fuzzy-number kind spaces. The theory we develop is used
to solve \textquotedblleft classical\textquotedblright\ fuzzy systems of
differential equations, including, for example, the fuzzy Cauchy problem and
the fuzzy wave equation. These tools allow us to obtain \emph{explicit}
solutions to fuzzy initial value problems which bear explicit formulas
similar to the crisp case, with some additional fuzzy terms which in the
crisp case disappear. The semigroup method displays a clear advantage over
other methods available in the literature (i.e., the level set method, the
differential inclusions method and other "fuzzification" methods of the
real-valued solution) in the sense that the solutions can be easily
constructed, and that the method can be applied to a larger class of fuzzy
differential equations that can be transformed into an abstract Cauchy
problem.
\end{abstract}


Subject Classification 2000 (AMS): 26E50, 46S40, 47N20, 47S40.

\section{Introduction}

\quad \quad Many physical processes in science are described by models using
differential equations. It is a well-known fact that the theory of
semigroups of operators on Banach spaces represents a powerful tool for
solving many classical differential equations (see e.g. $\left[ 9-10\right] $%
). Nowadays, differential equations are classified by the different
approaches used in order to deal with them: (i) deterministic (ii) random or
stochastic and (iii) fuzzy-like. The class (i) is by far the most studied.
The second class is usually employed when a deterministic description of the
model is not strictly justifiable. However, the first or second descriptions
often represent an idealization of real-world situations when imprecision or
lack of information in the modelling of the physical process may in fact
play a significant role. Indeed, it was this uncertainty or vagueness of the
modelling in certain physical processes that suggested the introduction of
the so-called fuzzy differential equations [3-4], whose solutions represent
fuzzy-number-valued functions.

The goal of this paper is to develop, in Section $3,$ the theory of the
semigroups of operators on spaces of fuzzy-number-valued functions. In
Section 4, our goal is to apply this theory to several classes of fuzzy
differential equations. The main difficulty of dealing with fuzzy
differential equations is the fact that the spaces of fuzzy numbers and
fuzzy-number-valued functions are not linear spaces. In particular, they are
not groups with respect to addition and the scalar multiplication is not, in
general, distributive with respect to usual addition of scalars. However,
they are complete metric spaces with their metrics having nice properties.
These features allow us to develop a consistent operator theory. To give a
better idea of what we are aiming for, let us recall that, as in the
classical cases (i) and (ii), many inhomogeneous fuzzy differential
equations (iii) can be recasted into the form of an abstract Cauchy problem
of the form%
\begin{align*}
\frac{du}{dt}\left( t\right) & =A[u\left( t\right) ]\oplus g\left( t\right) ,%
\text{ }t\in I, \\
u\left( t_{0}\right) & =u_{0},
\end{align*}%
where $I$ is an interval, $u:I\rightarrow X,$ $u_{0}\in X$, $A:X\rightarrow X
$ is a linear bounded operator on $X,$ and $g:I\rightarrow X$ is continuous,
where $(X,\oplus ,\odot ,d)$ is a complete metric space. Our goal is to show
that, as in the classical case, the following variation of parameter formula%
\begin{equation*}
u(t)=T(t)(u_{0})\oplus \int_{0}^{t}T(t-s)g(s)ds
\end{equation*}%
always furnishes a solution to the fuzzy Cauchy problem. Furthermore, we
also provide the mean for constructing $T(t)(u_{0})$ as follows:%
\begin{equation*}
T\left( t\right) (u_{0}):=e^{t\odot A}(u_{0})=\sum\limits_{k=0}^{\infty }%
\frac{t^{k}}{k!}A^{k}(u_{0}),
\end{equation*}%
where $\sum $ is the sum with respect to $\oplus .$ Finally, we apply this
formula on some specific examples and show that it always provides the
correct answer in the sense that the solutions we construct in the fuzzy
case contain additional terms which in the crisp case completely disappear.
It is worth mentioning that one can easily formulate the classical Liapunov
stability theorem for $\left\{ T\left( t\right) \right\} _{t\geq 0}$ in this
form. While the above formula provides a constructive mean to obtain
solutions, the following important questions still remain open: in the
generic case is the solution of the above Cauchy problem unique, and if so,
does it depend continuously on the initial data? 

\section{Preliminaries}

\quad \quad In this section we present the main properties of the space of
fuzzy numbers and of some other spaces based on it and which have similar
properties. Given a set $X\neq \emptyset ,$ a fuzzy subset of $X$ is a
mapping $u:X\rightarrow \left[ 0,1\right] $ and obviously any classical
subset of $X$ can be identified as a fuzzy subset of $X$ defined by $\chi
_{A}:X\rightarrow \left[ 0,1\right] ,\chi _{A}\left( x\right) =1$ if $x\in
A,\chi _{A}\left( x\right) =0,$ if $x\in X\setminus A.$ If $u:X\rightarrow %
\left[ 0,1\right] $ is a fuzzy subset of $X,$ then for $x\in X,$ $u\left(
x\right) $ is called the membership degree of $x$ to $u$ (see e.g. $\left[ 21%
\right] $).

\textbf{Definition 2.1}(see e.g. $\left[ 8\right] $). The space of fuzzy
numbers denoted by $\mathbf{R}_{\mathcal{F}}$ is defined as the class of
fuzzy subsets of the real axis $\mathbf{R},$ i.e., of $u:\mathbf{R}%
\rightarrow \left[ 0,1\right] ,$ having the following four properties:

$\left( i\right) $ $\forall u\in \mathbf{R}_{\mathcal{F}},$ $u$ is normal,
i.e., $\exists x_{u}\in \mathbf{R}$ with $u\left( x_{u}\right) =1;$

$\left( ii\right) $ $\forall u\in \mathbf{R}_{\mathcal{F}},$ $u$ is a convex
fuzzy set, i.e., 
\begin{equation*}
u\left( tx+\left( 1-t\right) y\right) \geq \min \left\{ u\left( x\right)
,u\left( y\right) \right\} ,\forall t\in \left[ 0,1\right] ,x,y\in \mathbf{R}%
;
\end{equation*}

$\left( iii\right) $ $\forall u\in \mathbf{R}_{\mathcal{F}},$ $u$ is
upper-semi-continuous on $\mathbf{R};$

$\left( iv\right) $ $\overline{\left\{ x\in \mathbf{R}\text{;}u\left(
x\right) >0\right\} }$ is compact, where $\overline{M}$ denotes the closure
of $M.$

\textbf{Remarks.} 1) Clearly, $\mathbf{R}\subset \mathbf{R}_{\mathcal{F}}$
because any real number $x_{0}\in \mathbf{R}$ can be identified with $\chi
_{\left\{ x_{0}\right\} },$ which satisfies the properties $\left( i\right)
-\left( iv\right) $ as in Definition $2.1.$

2) For $0<r\leq 1$ and $u\in \mathbf{R}_{\mathcal{F}},$ let us denote by $%
\left[ u\right] ^{r}=\left\{ x\in \mathbf{R};u\left( x\right) \geq r\right\} 
$ and $\left[ u\right] ^{0}=\overline{\left\{ x\in \mathbf{R};u\left(
x\right) >0\right\} },$ the so called level sets of $u.$ Then it is an
immediate consequence of $\left( i\right) -\left( iv\right) $ that $\left[ u%
\right] ^{r}$ represents a bounded closed (i.e., compact) subinterval of $%
\mathbf{R},$ denoted by $\left[ u\right] ^{r}=\left[ u_{-}\left( r\right)
,u_{+}\left( r\right) \right] ,$ where $u_{-}\left( r\right) \leq
u_{+}\left( r\right) $ for all $r\in \left[ 0,1\right] .$

\textbf{Definition 2.2} (see e.g. $\left[ 8\right] $). The addition and the
scalar product in $\mathbf{R}_{\mathcal{F}}$ are defined by $\oplus :\mathbf{%
R}_{\mathcal{F}}\times \mathbf{R}_{\mathcal{F}}\rightarrow \mathbf{R}_{%
\mathcal{F}},$%
\begin{equation*}
\left( u\oplus v\right) \left( x\right) =\sup_{y+z=x}\min \left\{ u\left(
y\right) ,v\left( z\right) \right\}
\end{equation*}%
and, by $\odot :\mathbf{R}\times \mathbf{R}_{\mathcal{F}}\rightarrow \mathbf{%
R}_{\mathcal{F}},$%
\begin{equation*}
\left( \lambda \odot v\right) \left( x\right) =\left\{ 
\begin{array}{l}
u\left( \frac{x}{\lambda }\right) \\ 
\widetilde{0}%
\end{array}%
\begin{array}{l}
\text{if }\lambda \neq 0, \\ 
\text{if }\lambda =0,%
\end{array}%
\right.
\end{equation*}%
where $\widetilde{0}:\mathbf{R}\rightarrow \left[ 0,1\right] $ is $%
\widetilde{0}=\chi _{\left\{ 0\right\} }.$

Also, we can write%
\begin{equation*}
\left[ u\oplus v\right] ^{r}=\left[ u\right] ^{r}+\left[ v\right] ^{r},\left[
\lambda \odot v\right] ^{r}=\lambda \left[ v\right] ^{r},
\end{equation*}%
for all $r\in \left[ 0,1\right] ,$ where $\left[ u\right] ^{r}+\left[ v%
\right] ^{r}$ means the usual sum of two intervals (as subsets of $\mathbf{R}
$) and $\lambda \left[ v\right] ^{r}$ means the usual product between a real
scalar and a subset of $\mathbf{R}$. If we define $D:\mathbf{R}_{\mathcal{F}%
}\times \mathbf{R}_{\mathcal{F}}\rightarrow \mathbf{R}_{+}\cup \left\{
0\right\} $ by 
\begin{equation*}
D\left( u,v\right) =\sup_{r\in \left[ 0,1\right] }\max \left\{ \left\vert
u_{-}\left( r\right) -v_{-}\left( r\right) \right\vert ,\left\vert
u_{+}\left( r\right) -v_{+}\left( r\right) \right\vert \right\} ,
\end{equation*}%
where $\left[ u\right] ^{r}=\left[ u_{-}\left( r\right) ,u_{+}\left(
r\right) \right] ,$ $\left[ v\right] ^{r}=\left[ v_{-}\left( r\right)
,v_{+}\left( r\right) \right] ,$ then we have the following.

\textbf{Theorem 2.3} ( see e.g. $\left[ 8\right] $). \textit{$\left( \mathbf{%
R}_{\mathcal{F}},D\right) $ is a complete metric space. In addition, $D$ has
the following three properties: }

\textit{$\left( i\right) $ $D\left( u\oplus w, v\oplus w\right) =D\left(
u,v\right) ,$ for all $u,v,w\in \mathbf{R}_{\mathcal{F}};$ }

\textit{$\left( ii\right) $ $D\left( k\odot u,k\odot v\right) =\left|
k\right| D\left( u,v\right) ,$ for all $u,v\in \mathbf{R}_{\mathcal{F}},k\in 
\mathbf{R};$ }

\textit{$\left( iii\right) $ $D\left( u\oplus v,w\oplus e\right) \leq
D\left( u,w\right) +D\left( v,e\right) ,$ for all $u,v,w,e\in \mathbf{R}_{%
\mathcal{F}}.$}

Also, the following result is known.

\textbf{Theorem 2.4} (see e.g. $\left[ 1\right] ,\left[ 8\right] $). \textit{%
$\left( i\right) $ $u\oplus v=v\oplus u,u\oplus \left( v\oplus w\right)
=\left( u\oplus v\right) \oplus w;$ }

\textit{$\left( ii\right) $ If we denote $\widetilde{0}=\chi _{\left\{
0\right\} },$ then $u\oplus \widetilde{0}=\widetilde{0}\oplus u=u,$ for any $%
u\in \mathbf{R}_{\mathcal{F}};$ }

\textit{$\left( iii\right) $ With respect to $\widetilde{0},$ none of $u\in 
\mathbf{R}_{\mathcal{F}}\setminus \mathbf{R}$ has an opposite member
(regarding $\oplus $) in $\mathbf{R}_{\mathcal{F}};$ }

\textit{$\left( iv\right) $ for any $a,b\in \mathbf{R}$ with $a,b\geq 0$ or $%
a,b\leq 0$ and any $u\in \mathbf{R}_{\mathcal{F}},$ we have 
\begin{equation*}
\left( a+b\right) \odot u=a\odot u\oplus b\odot u.
\end{equation*}
For general $a,b\in \mathbf{R},$ the above property does not hold. }

\textit{$\left( v\right) $ $\lambda \odot \left( u\oplus v\right) =\lambda
\odot u\oplus \lambda \odot v,$ for all $\lambda \in \mathbf{R},u,v\in 
\mathbf{R}_{\mathcal{F}};$ }

\textit{$\left( vi\right) $ $\lambda \odot \left( \mu \odot u\right) =\left(
\lambda \mu \right) \odot u,$ for all $\lambda ,\mu \in \mathbf{R},u\in 
\mathbf{R}_{\mathcal{F}};$ }

\textit{$\left( vii\right) $ If we denote $\left\| u\right\| _{\mathcal{F}%
}=D\left( u,\widetilde{0}\right) ,u\in \mathbf{R}_{\mathcal{F}},$ then $%
\left\| u\right\| _{\mathcal{F}}$ has the properties of an usual norm on $%
\mathbf{R}_{\mathcal{F}},$ i.e. $\left\| u\right\| _{\mathcal{F}}=0$ iff $u=%
\widetilde{0},\left\| \lambda \odot u\right\| _{\mathcal{F}}=\left| \lambda
\right| \left\| u\right\| _{\mathcal{F}},\left\| u+v\right\| _{\mathcal{F}%
}\leq \left\| u\right\| _{\mathcal{F}}+\left\| v\right\| _{\mathcal{F}%
},\left| \left\| u\right\| _{\mathcal{F}}-\left\| v\right\| _{\mathcal{F}%
}\right| \leq D\left( u,v\right) ;$ }

\textit{$\left( viii\right) D\left( \alpha \odot u,\beta \odot u\right)
=\left| \alpha -\beta \right| D\left( \widetilde{0},u\right) ,$ for all $%
\alpha ,\beta \geq 0,u\in \mathbf{R}_{\mathcal{F}}.$ If $\alpha ,\beta \leq
0 $ then the equality is also valid. If $\alpha $ and $\beta $ are of
opposite signs, then the equality is not valid.}

\textbf{Remarks.} 1) Theorem $2.4$ shows that $\left( \mathbf{R}_{\mathcal{F}%
},\oplus ,\odot \right) $ is not a linear space over $\mathbf{R}$ and,
consequently, $\left( \mathbf{R}_{\mathcal{F}},\left\Vert u\right\Vert _{%
\mathcal{F}}\right) $ cannot be a normed space.

2) On $\mathbf{R}_{\mathcal{F}},$ we can define the substraction $\ominus ,$
called the $H-$ difference (see $\left[ 12\right] $) as follows: $u\ominus v$
has sense if there exists $w\in \mathbf{R}_{\mathcal{F}}$ such that $%
u=v\oplus w.$ Clearly, $u\ominus v$ does not exist for all $u,v\in \mathbf{R}%
_{\mathcal{F}}$ (for example, $\widetilde{0}\ominus v$ does not exists if $%
v\neq \widetilde{0}$).

In what follows, we define some usual spaces of fuzzy-number-valued
functions, which have similar properties to $\left( \mathbf{R}_{\mathcal{F}%
},D\right) .$

Denote $C\left( \left[ a,b\right] ;\mathbf{R}_{\mathcal{F}}\right) =\left\{
f:\left[ a,b\right] \rightarrow \mathbf{R}_{\mathcal{F}};\text{ }f\text{ is
continuous on }\left[ a,b\right] \right\} ,$ endowed with the metric 
\begin{equation*}
D^{\ast }\left( f,g\right) =\sup \left\{ D\left( f\left( x\right) ,g\left(
x\right) \right) ;x\in \left[ a,b\right] \right\} .
\end{equation*}%
Since $\left( \mathbf{R}_{\mathcal{F}},D\right) $ is a complete metric
space, by standard technique (see e.g. $\left[ 20\right] $) we obtain that $%
\left( C\left( \left[ a,b\right] ;\mathbf{R}_{\mathcal{F}}\right) ,D^{\ast
}\right) $ is also complete metric space. Moreover, if we define 
\begin{equation*}
\left( f\oplus g\right) \left( x\right) =f\left( x\right) \oplus g\left(
x\right) ,\left( \lambda \odot f\right) \left( x\right) =\lambda \odot
f\left( x\right) ,
\end{equation*}%
(for simplicity, the addition and scalar multiplication in $C\left( \left[
a,b\right] ;\mathbf{R}_{\mathcal{F}}\right) $ are denoted as in $\mathbf{R}_{%
\mathcal{F}}$), also $\widetilde{0}:\left[ a,b\right] \rightarrow \mathbf{R}%
_{\mathcal{F}},\widetilde{0}\left( t\right) =\widetilde{0}_{\mathbf{R}_{%
\mathcal{F}}},$ for all $t\in \left[ a,b\right] ,$%
\begin{equation*}
\left\Vert f\right\Vert _{\mathcal{F}}=\sup \left\{ D\left( \widetilde{0}%
,f\left( x\right) \right) ;x\in \left[ a,b\right] \right\} ,
\end{equation*}%
then we easily obtain the following properties.

\textbf{Theorem 2.5} \textit{$\left( i\right) $ $f\oplus g=g\oplus f,\left(
f\oplus g\right) \oplus h=f\oplus \left( g\oplus h\right) ;$ }

\textit{$\left( ii\right) $ $f\oplus \widetilde{0}=\widetilde{0}\oplus f,$
for any $f\in C\left( \left[ a,b\right] ;\mathbf{R}_{\mathcal{F}}\right) ;$ }

\textit{$\left( iii\right) $ With respect to $\widetilde{0}$ in $C\left( %
\left[ a,b\right] ;\mathbf{R}_{\mathcal{F}}\right) ,$ any $f\in C\left( %
\left[ a,b\right] ;\mathbf{R}_{\mathcal{F}}\right) $ with $f\left( \left[ a,b%
\right] \right) \cap \mathbf{R}_{\mathcal{F}}\neq \emptyset $ has no
opposite member (regarding $\oplus $) in $C\left( \left[ a,b\right] ;\mathbf{%
R}_{\mathcal{F}}\right) ;$ }

\textit{$\left( iv\right) $ for all $\lambda ,\mu \in \mathbf{R}$ with $%
\lambda ,\mu \geq 0$ or $\lambda ,\mu \leq 0$ and for any $f\in C\left( %
\left[ a,b\right] ;\mathbf{R}_{\mathcal{F}}\right) ,$%
\begin{equation*}
\left( \lambda +\mu \right) \odot f=\left( \lambda \odot f\right) \oplus
\left( \mu \odot f\right) ;
\end{equation*}
For general $\lambda ,\mu \in \mathbf{R},$ this property does not hold. }

\textit{$\left( v\right) $ $\lambda \odot \left( f\oplus g\right) =\lambda
\odot f\oplus \lambda \odot g,\lambda \odot \left( \mu \odot f\right)
=\left( \lambda \mu \right) \odot f,$ for any $f,g\in C\left( \left[ a,b%
\right] ;\mathbf{R}_{\mathcal{F}}\right) ,\lambda ,\mu \in \mathbf{R};$ }

\textit{$\left( vi\right) $ $\left\| f\right\| _{\mathcal{F}}=0$ iff $f=%
\widetilde{0},$ $\left\| \lambda \odot f\right\| _{\mathcal{F}}=\left|
\lambda \right| \left\| f\right\| _{\mathcal{F}},\left\| f\oplus g\right\| _{%
\mathcal{F}}\leq \left\| f\right\| _{\mathcal{F}}+\left\| g\right\| _{%
\mathcal{F}},\left| \left\| f\right\| _{\mathcal{F}}-\left\| g\right\| _{%
\mathcal{F}}\right| \leq D^{*}\left( f,g\right) ,$ for any $f,g\in C\left( %
\left[ a,b\right] ;\mathbf{R}_{\mathcal{F}}\right) ,\lambda \in \mathbf{R};$ 
}

\textit{$\left( vii\right) $ $D^{*}\left( \lambda \odot f,\mu \odot f\right)
=\left| \lambda -\mu \right| D^{*}\left( \widetilde{0},f\right) ,$for any $%
f\in C\left( \left[ a,b\right] ;\mathbf{R}_{\mathcal{F}}\right) ,\lambda \mu
\geq 0;$ }

\textit{$\left( viii\right) $%
\begin{equation*}
D^{*}\left( f\oplus h,g\oplus h\right) =D^{*}\left( f,g\right) ,
\end{equation*}
\begin{equation*}
D^{*}\left( \lambda \odot f,\lambda \odot g\right) =\left| \lambda \right|
D^{*}\left( f,g\right) ,
\end{equation*}
\begin{equation*}
D^{*}\left( f\oplus g,h\oplus e\right) \leq D^{*}\left( f,h\right)
+D^{*}\left( g,e\right) ,
\end{equation*}
for any $f,g,h,e\in C\left( \left[ a,b\right] ;\mathbf{R}_{\mathcal{F}%
}\right) ,\lambda \in \mathbf{R}.$}

\textbf{Remark.} It is easy to show that if $f,g\in C\left( \left[ a,b\right]
;\mathbf{R}_{\mathcal{F}}\right) ,$ then $F:\left[ a,b\right] \rightarrow 
\mathbf{R},$ defined by $F\left( x\right) =D\left( f\left( x\right) ,g\left(
x\right) \right) $ is continuous on $\left[ a,b\right] .$

Now, for $1\leq p<\infty $ and a strongly measurable function $f$ on $\left[
a,b\right] ,$ let us define 
\begin{equation*}
L^{p}\left( \left[ a,b\right] ;\mathbf{R}_{\mathcal{F}}\right) =\text{ }%
\left\{ f:\left( L\right) \int\limits_{a}^{b}\left( D\left( \widetilde{0}%
,f\left( x\right) \right) \right) ^{p}dx<+\infty \right\} ,
\end{equation*}%
where according to e.g. $\left[ 7\right] $, $f$ is called strongly
measurable if, for each $x\in \left[ a,b\right] ,$ $f_{-}\left( x\right)
\left( r\right) $ and $f_{+}\left( x\right) \left( r\right) $ are Lebesgue
measurable as functions of $r\in \left[ 0,1\right] $ (here again, $\left[
f\left( x\right) \right] ^{r}=\left[ f_{-}\left( x\right) \left( r\right)
,f_{+}\left( x\right) \left( r\right) \right] $ denotes the $r$-level set of 
$f\left( x\right) \in \mathbf{R}_{\mathcal{F}}$). The following result shows
that $L^{p}\left( \left[ a,b\right] ;\mathbf{R}_{\mathcal{F}}\right) $ is
well defined.

\textbf{Theorem 2.6} \textit{$\left( i\right) $ If $f:\left[ a,b\right]
\rightarrow \mathbf{R}_{\mathcal{F}}$ is strongly measurable then $F:\left[
a,b\right] \rightarrow \mathbf{R}_{+}$ defined by $F\left( x\right) =D\left( 
\widetilde{0},f\left( x\right) \right) $ is Lebesgue measurable on $\left[
a,b\right] ;$ }

\textit{$\left( ii\right) $ For any $f,g\in L^{p}\left( \left[ a,b\right] ;%
\mathbf{R}_{\mathcal{F}}\right) ,F\left( x\right) =D\left( f\left( x\right)
,g\left( x\right) \right) $ is Lebesgue measurable and $L^{p}$-integrable on 
$\left[ a,b\right] .$ Moreover, if we define 
\begin{equation*}
D_{p}\left( f,g\right) =\left\{ \left( L\right) \int\limits_{a}^{b}\left[
D\left( \widetilde{0},f\left( x\right) \right) \right] ^{p}dx\right\} ^{%
\frac{1}{p}},
\end{equation*}
then $\left( L^{p}\left( \left[ a,b\right] ;\mathbf{R}_{\mathcal{F}}\right)
,D_{p}\right) $ is a complete metric space (in $L^{p}\left( \left[ a,b\right]
;\mathbf{R}_{\mathcal{F}}\right) ,$ $f=g$ means $f\left( x\right) =g\left(
x\right) ,$ a.e. $x\in \left[ a,b\right] $) and, in addition, $D_{p} $
satisfies the following properties: 
\begin{equation*}
D_{p}\left( f\oplus h,g\oplus h\right) =D_{p}\left( f,g\right) ,
\end{equation*}
\begin{equation*}
D_{p}\left( \lambda \odot f,\lambda \odot g\right) =\left| \lambda \right|
D_{p}\left( f,g\right) ,
\end{equation*}
\begin{equation*}
D_{p}\left( f\oplus g,h\oplus e\right) \leq D_{p}\left( f,h\right)
+D_{p}\left( g,e\right) ,
\end{equation*}
for any $f,g,h,e\in L^{p}\left( \left[ a,b\right] ;\mathbf{R}_{\mathcal{F}%
}\right) ,\lambda \in \mathbf{R}.$}

\textbf{Proof.} $\left( i\right) $ By definition, we have 
\begin{equation*}
D\left( \widetilde{0},f\left( x\right) \right) =\sup_{r\in \left[ 0,1\right]
}\max \left\{ \left\vert f_{-}\left( x\right) \left( r\right) \right\vert
,\left\vert f_{+}\left( x\right) \left( r\right) \right\vert \right\} ,
\end{equation*}%
for all $x\in \left[ a,b\right] .$ Because $f_{-}\left( x\right) \left(
r\right) $ and $f_{+}\left( x\right) \left( r\right) $ are left continuous
at each $r\in \left[ 0,1\right] $ and right continuous at $r=0$ (see e.g. $%
\left[ 11\right] $), denoting by $\mathbb{Q}$, the set of all rational
numbers, we easily obtain 
\begin{equation*}
D\left( \widetilde{0},f\left( x\right) \right) =\sup_{r\in \left[ 0,1\right]
\cap \mathbf{Q}}\max \left\{ \left\vert f_{-}\left( x\right) \left( r\right)
\right\vert ,\left\vert f_{+}\left( x\right) \left( r\right) \right\vert
\right\} ,
\end{equation*}%
which implies that $D\left( \widetilde{0},f\left( x\right) \right) $ is
Lebesgue measurable on $\left[ a,b\right] $ (as supremum of a countable set
of Lebesgue measurable functions).

$\left( ii\right) $ We have 
\begin{align*}
F\left( x\right) & =D\left( f\left( x\right) ,g\left( x\right) \right) \\
& =\sup_{r\in \left[ 0,1\right] }\max \left\{ \left\vert f_{-}\left(
x\right) \left( r\right) -g_{-}\left( x\right) \left( r\right) \right\vert
,\left\vert f_{+}\left( x\right) \left( r\right) -g_{+}\left( x\right)
\left( r\right) \right\vert \right\} ,
\end{align*}%
and reasoning exactly as for the above point $\left( i\right) $, it follows
that $F$ is Lebesgue measurable on $\left[ a,b\right] .$ Moreover, since 
\begin{equation*}
D\left( f\left( x\right) ,g\left( x\right) \right) \leq D\left( f\left(
x\right) ,\widetilde{0}\right) +D\left( \widetilde{0},g\left( x\right)
\right) ,
\end{equation*}%
we immediately obtain%
\begin{equation*}
\left( L\right) \int\limits_{a}^{b}\left[ D\left( f\left( x\right) ,g\left(
x\right) \right) \right] ^{p}dx<+\infty .
\end{equation*}%
Then by the properties of $D,$ it easily follows that $D_{p}$ is a metric on 
$L^{p}\left( \left[ a,b\right] ;\mathbf{R}_{\mathcal{F}}\right) $ satisfying
the properties listed above at $\left( ii\right) .$

It remains to prove that $\left( L^{p}\left( \left[ a,b\right] ;\mathbf{R}_{%
\mathcal{F}}\right) ,D_{p}\right) $ is complete. Let $f_{n}\in L^{p}\left( %
\left[ a,b\right] ;\mathbf{R}_{\mathcal{F}}\right) ,n\in \mathbb{N}$ be a
Cauchy sequence with respect to the metric $D_{p}.$ By standard technique
(see e.g. $\left[ 20\right] $), we obtain a subsequence $\left(
f_{n_{j}}\right) _{j\in \mathbb{N}}$ such that $f_{n_{j}}\left( x\right)
,j\in \mathbb{N},$ a.e. $x\in \left[ a,b\right] $ is a Cauchy sequence in
the complete metric space $\left( \mathbf{R}_{\mathcal{F}},D\right) .$
Therefore, there exists $f\left( x\right) $ (a.e. $x\in \left[ a,b\right] $)
such that 
\begin{equation*}
\lim_{j\rightarrow +\infty }D\left( f_{n_{j}}\left( x\right) ,f\left(
x\right) \right) =0.
\end{equation*}%
Then by e.g. $[7,$Proposition $6.1.4]$, it follows that $f$ is strongly
measurable (on $\left[ a,b\right] $). The rest of the proof follows standard
technique as in e.g. $\left[ 20\right] .$ This proves the theorem. $\Box $

In what follows, we introduce the following definition.

\textbf{Definition 2.7} (i) (see e.g. $\left[ 7\right] $). Let $f:\left[ a,b%
\right] \rightarrow \mathbf{R}_{\mathcal{F}}$ and $x_{0}\in \left(
a,b\right) .$ We say that $f$ is (Hukuhara) differentiable at $x_{0},$ if
there exists $h_{0}>0$ such that for all $0<h\leq h_{0},$ there exists $%
f\left( x_{0}+h\right) \ominus f\left( x_{0}\right) ,f\left( x_{0}\right)
\ominus f\left( x_{0}-h\right) $ and $f^{^{\prime }}\left( x_{0}\right) \in 
\mathbf{R}_{\mathcal{F}}$ with the property: 
\begin{align*}
\lim_{h\searrow 0}D\left[ \frac{1}{h}\odot \left( f\left( x_{0}+h\right)
\ominus f\left( x_{0}\right) \right) ,f^{^{\prime }}\left( x_{0}\right) %
\right] & =0, \\
\lim_{h\searrow 0}D\left[ \frac{1}{h}\odot \left( f\left( x_{0}\right)
\ominus f\left( x_{0}-h\right) \right) ,f^{^{\prime }}\left( x_{0}\right) %
\right] & =0.
\end{align*}

(ii) (see [2]). Let $f:(a,b)\rightarrow \mathbf{R}_{\mathcal{F}}$ and $t\in
(a,b).$ We say that $f$ is generalized differentiable at $t,$ if:

1) There exist $f(t+h)\ominus f(t)$, $f(t)\ominus f(t-h)$, for all $h>0$
sufficiently small and there exist 
\begin{equation*}
\lim_{h\searrow 0}\frac{f(t+h)\ominus f(t)}{h}=\lim_{h\searrow 0}\frac{%
f(t)\ominus f(t-h)}{h}=f^{\prime }(t)\in \mathbf{R}_{\mathcal{F}}
\end{equation*}
or

2) There exist $f(t)\ominus f(t+h)$, $f(t-h)\ominus f(t)$, for all $h>0$
sufficiently small and there exist 
\begin{equation*}
\lim_{h\searrow 0}\frac{f(t)\ominus f(t+h)}{-h}=\lim_{h\searrow 0}\frac{%
f(t-h)\ominus f(t)}{-h}=f^{\prime }(t)\in \mathbf{R}_{\mathcal{F}}
\end{equation*}
or

3) There exist $f(t+h)\ominus f(t)$, $f(t-h)\ominus f(t)$, for all $h>0$
sufficiently small and there exist 
\begin{equation*}
\lim_{h\searrow 0}\frac{f(t+h)\ominus f(t)}{h}=\lim_{h\searrow 0}\frac{%
f(t-h)\ominus f(t)}{-h}=f^{\prime }(t)\in \mathbf{R}_{\mathcal{F}}
\end{equation*}
or

4) There exist $f(t)\ominus f(t-h)$, $f(t)\ominus f(t+h)$, for all $h>0$
sufficiently small and there exist 
\begin{equation*}
\lim_{h\searrow 0}\frac{f(t)\ominus f(t-h)}{h}=\lim_{h\searrow 0}\frac{%
f(t)\ominus f(t+h)}{-h}=f^{\prime }(t)\in \mathbf{R}_{\mathcal{F}}
\end{equation*}
(Here all the limits are considered in the metric $D$ and $h$ at
denominators means $\frac{1}{h}\odot $). Obviously the Hukuhara
differentiability implies the generalized differentiability but the converse
is not valid.

If $f^{^{\prime }}\left( x_{0}\right) $ exists for all $x_{0}\in \left[ a,b%
\right] ,$ by iteration we can define%
\begin{equation*}
f^{^{\prime \prime }}\left( x_{0}\right) =\left( f^{^{\prime }}\right)
^{^{\prime }}\left( x_{0}\right) ,x_{0}\in \left[ a,b\right] .
\end{equation*}%
For $p\in \mathbb{N}$, let us denote by 
\begin{equation*}
C^{p}\left( \left[ a,b\right] ;\mathbf{R}_{\mathcal{F}}\right) =\left\{ f%
\text{ }:\left[ a,b\right] \rightarrow \mathbf{R}_{\mathcal{F}}\text{; there
exists }f^{\left( p\right) }\text{ continuous on }\left[ a,b\right] \text{ }%
\right\}
\end{equation*}%
(here the derivatives $f^{\prime },f",...,f^{(p)}$ are all considered of the
same kind as in Definition 2.7) and for $f,g\in C^{p}\left( \left[ a,b\right]
;\mathbf{R}_{\mathcal{F}}\right) ,$ let us define 
\begin{equation*}
D_{p}^{\ast }\left( f,g\right) =\sum\limits_{i=0}^{p}D^{\ast }\left(
f^{\left( i\right) },g^{\left( i\right) }\right) ,
\end{equation*}%
where%
\begin{equation*}
D^{\ast }\left( f^{\left( i\right) },g^{\left( i\right) }\right) =\sup
\left\{ D\left( f^{\left( i\right) }\left( x\right) ,g^{\left( i\right)
}\left( x\right) \right) ;x\in \left[ a,b\right] \right\} .
\end{equation*}%
Then we obtain that $D_{p}^{\ast }\left( f,g\right) $ is a metric and by
standard technique (as in the case of the real-valued functions), we obtain
that $\left( C^{p}\left( \left[ a,b\right] ;\mathbf{R}_{\mathcal{F}}\right)
,D_{p}^{\ast }\right) $ is a complete metric space. In addition, $%
D_{p}^{\ast }$ shares the same properties with $D^{\ast }$ (see, Theorem $%
2.5,$ $\left( vii\right) $ and $\left( viii\right) $).

Other spaces with properties similar to those of $\left( \mathbf{R}_{%
\mathcal{F}},D\right) $ can be constructed as follows. For $p\geq 1,$ let us
define 
\begin{equation*}
l_{\mathbf{R}_{\mathcal{F}}}^{p}=\left\{ x=\left( x_{n}\right) _{n};x_{n}\in 
\mathbf{R}_{\mathcal{F}},\forall n\in \mathbb{N}\text{ and }%
\sum\limits_{n=1}^{\infty }\left\Vert x_{n}\right\Vert _{\mathbf{R}_{%
\mathcal{F}}}^{p}<+\infty \right\} ,
\end{equation*}%
endowed with the metric 
\begin{equation*}
\rho _{p}\left( x,y\right) =\left\{ \sum\limits_{n=1}^{\infty }[D\left(
x_{n},y_{n}\right) ]^{p}\right\} ^{\frac{1}{p}},\forall x=\left(
x_{n}\right) _{n},y=\left( y_{n}\right) _{n}\in l_{\mathbf{R}_{\mathcal{F}%
}}^{p}.
\end{equation*}%
By the inequality,%
\begin{align*}
D\left( x_{n},y_{n}\right) & =D\left( x_{n}\oplus \widetilde{0},\widetilde{0}%
\oplus y_{n}\right) \\
& \leq D\left( x_{n},\widetilde{0}\right) +D\left( \widetilde{0},y_{n}\right)
\\
& =\left\Vert x_{n}\right\Vert _{\mathcal{F}}+\left\Vert y_{n}\right\Vert _{%
\mathcal{F}},
\end{align*}%
we easily get (by Minkowski's inequality if $p>1$) $\rho _{p}\left(
x,y\right) <+\infty .$ Also, it easily follows that $\rho _{p}\left(
x,y\right) $ is a metric with similar properties to $D$ (see Theorem $2.3$
and Theorem 2.4, (vii), (viii)). Since again $\left( \mathbf{R}_{\mathcal{F}%
},D\right) $ is a complete metric space, by the standard technique, we
easily get that $\left( l_{\mathbf{R}_{\mathcal{F}}}^{p},\rho _{p}\right) $
is also a complete metric space. Next, let us denote%
\begin{equation*}
m_{\mathbf{R}_{\mathcal{F}}}=\left\{ x=\left( x_{n}\right) _{n};x_{n}\in 
\mathbf{R}_{\mathcal{F}},\forall n\in \mathbb{N}\text{ and }\exists M>0\text{
such that }\left\Vert x_{n}\right\Vert _{\mathcal{F}}\leq M,\forall n\in 
\mathbb{N}\right\} ,
\end{equation*}%
endowed with the metric 
\begin{equation*}
\mu \left( x,y\right) =\sup \left\{ D\left( x_{n},y_{n}\right) ,\forall n\in 
\mathbb{N}\right\} .
\end{equation*}%
We easily get that $\left( m_{\mathbf{R}_{\mathcal{F}}},\mu \right) $ is a
complete metric space and, in addition, $\mu $ has similar properties to $D$
(see Theorem $2.3$ and Theorem 2.4, (vii), (viii)). Similarly, if we denote 
\begin{equation*}
c_{\mathbf{R}_{\mathcal{F}}}=\left\{ x=\left( x_{n}\right) _{n};x_{n}\in 
\mathbf{R}_{\mathcal{F}},\forall n\in \mathbb{N}\text{ and }\exists a\in 
\mathbf{R}_{\mathcal{F}}\text{ such that }D\left( x_{n},a\right) \overset{%
n\rightarrow \infty }{\longrightarrow }0\right\}
\end{equation*}%
and 
\begin{equation*}
c_{\mathbf{R}_{\mathcal{F}}}^{\widetilde{0}}=\left\{ x=\left( x_{n}\right)
_{n};x_{n}\in \mathbf{R}_{\mathcal{F}},\forall n\in \mathbb{N}\text{, such
that }D\left( x_{n},\widetilde{0}\right) \overset{n\rightarrow \infty }{%
\longrightarrow }0\right\} ,
\end{equation*}%
since $\left( \mathbf{R}_{\mathcal{F}},D\right) $ is complete, by standard
technique, it follows that $\left( c_{\mathbf{R}_{\mathcal{F}}},\mu \right) $
and $(c_{\mathbf{R}_{\mathcal{F}}}^{\widetilde{0}},\mu )$ are complete
metric spaces. Note that any finite Cartesian product of the spaces in this
section, endowed with the \textquotedblright box\textquotedblright\ metric
(i.e., $d=\underset{i}{\max }\rho _{i}$) is a complete metric space and the
\textquotedblright Cartesian\textquotedblright\ metric $d$ shares the same
properties with the metric $D$ on $\mathbf{R}_{\mathcal{F}}$. Clearly, the
concept of differentiability in Definition 2.7 remains valid if $\mathbf{R}_{%
\mathcal{F}}$ is replaced by any of the spaces $(X,\oplus ,\odot ,d)$
introduced above.

\section{Operator theory and semigroups of operators}

\quad \quad In this section, we consider elements of operator theory on the
complete metric spaces introduced in Section $2$. Everywhere in this section
one can think of $\left( X,\oplus ,\odot ,d\right) $ as any of the spaces
introduced in the previous section.

\textbf{Definition 3.1.} $A:X\mathbb{\rightarrow }\mathbf{R}$ is a linear
functional if 
\begin{equation*}
\left\{ 
\begin{array}{l}
A\left( x\oplus y\right) =A\left( x\right) +A\left( y\right) , \\ 
A\left( \lambda \odot x\right) =\lambda A\left( x\right) ,%
\end{array}%
\right.
\end{equation*}%
and, also, $A:X\mathbb{\rightarrow }X$ is a linear operator if 
\begin{equation*}
\left\{ 
\begin{array}{l}
A\left( x\oplus y\right) =A\left( x\right) \oplus A\left( y\right) , \\ 
A\left( \lambda \odot x\right) =\lambda \odot A\left( x\right) ,%
\end{array}%
\right.
\end{equation*}%
for all $x,y\in X$,$\lambda \in \mathbf{R}$.

\textbf{Remark.} If $A:X\mathbb{\rightarrow }\mathbf{R}$ or $A:X\mathbb{%
\rightarrow }X$ is linear and continuous at $\widetilde{0}\in X$, then the
latter does not imply the continuity of $A$ at each $x\in X$, because in
general, we cannot write $x_{0}=(x_{0}\ominus x)\oplus x$.

However, we can prove the following theorem.

\textbf{Theorem 3.2} \textit{$\left( i\right) $ If $A:X\mathbb{\rightarrow }%
\mathbf{R}$ is linear, then it is continuous at $\widetilde{0}\in X$ if and
only if there exists $M>0$ such that 
\begin{equation*}
\left\vert A\left( x\right) \right\vert \leq M\left\Vert x\right\Vert _{%
\mathcal{F}},\forall x\in X\text{,}
\end{equation*}%
where $\left\Vert x\right\Vert _{\mathcal{F}}=d\left( \widetilde{0},x\right)
.$ }

\textit{$\left( ii\right) $ If $A:X\mathbb{\rightarrow }X$ is linear, then
it is continuous at $\widetilde{0}\in X$ if and only if there exists $M>0$
such that 
\begin{equation*}
\left\Vert A\left( x\right) \right\Vert _{\mathcal{F}}\leq M\left\Vert
x\right\Vert _{\mathcal{F}},\forall x\in X.
\end{equation*}%
}

\textbf{Proof.} The proof is similar to that in the case of the usual Banach
spaces (see e.g. $\left[ 20\right] $) by taking into account the properties
of $\left( X,\oplus ,\odot ,d\right) .$ For example, let us prove $\left(
ii\right) .$ By the continuity of $A:X\mathbb{\rightarrow }X$ at $\widetilde{%
0},$ for $B\left( \widetilde{0},1\right) =\left\{ y\in X\text{; }d\left( 
\widetilde{0},y\right) <1\right\} ,$ there exists $V\in \mathcal{V}\left( 
\widetilde{0}\right) $ such that $A\left( V\right) \subset B\left( 
\widetilde{0},1\right) .$ Also, since $V\in \mathcal{V}\left( \widetilde{0}%
\right) ,$ $\exists r>0$ such that $B\left( \widetilde{0},r\right) =\left\{
y\in X\text{; }d\left( \widetilde{0},y\right) <r\right\} \subset V.$ It
follows $A\left( B\left( \widetilde{0},r\right) \right) \subset B\left( 
\widetilde{0},1\right) .$ Let $M=\frac{2}{r}$ and let $x\in X$ be arbitrary.
If $x=\widetilde{0},$ then by additivity, we get $A\left( \widetilde{0}%
\right) =\widetilde{0},$ and the inequality is trivial. Now let $x\neq 
\widetilde{0}.$ We get $\left\Vert x\right\Vert _{\mathcal{F}}=d\left( 
\widetilde{0},x\right) \neq 0$ and let us define $\alpha =\frac{r}{%
2\left\Vert x\right\Vert _{\mathcal{F}}}>0.$ We have $\left\Vert \alpha
\odot x\right\Vert _{\mathcal{F}}=\alpha \left\Vert x\right\Vert _{\mathcal{F%
}}=\frac{r}{2}<r,$ i.e. $\alpha \odot x\in B\left( \widetilde{0},r\right) .$
It follows $A\left( \alpha \odot x\right) \in B\left( \widetilde{0},1\right)
,$ i.e., $\left\Vert A\left( \alpha \odot x\right) \right\Vert _{\mathcal{F}%
}=d\left( \widetilde{0},A\left( \alpha \odot x\right) \right) <1,$ which
implies $\alpha \left\Vert A\left( x\right) \right\Vert _{\mathcal{F}}<1,$
i.e. $\left\Vert A\left( x\right) \right\Vert _{\mathcal{F}}<\frac{1}{\alpha 
}=\frac{2}{r}\left\Vert x\right\Vert _{\mathcal{F}}$ $=M\left\Vert
x\right\Vert _{\mathcal{F}}$. The converse in $\left( ii\right) $ is
immediate. $\Box $

Now, for $A:X\mathbb{\rightarrow }\mathbf{R}$ linear and continuous at $%
\widetilde{0},$ let us denote 
\begin{equation*}
\mathcal{M}_{A}:=\left\{ M>0;\left\vert A\left( x\right) \right\vert \leq
M\left\Vert x\right\Vert _{\mathcal{F}},\forall x\in X\right\} ,
\end{equation*}%
and for $A:X\mathbb{\rightarrow }X$ linear and continuous at $\widetilde{0},$
let us denote 
\begin{equation*}
\mathcal{M}_{A}:=\left\{ M>0;\left\Vert A\left( x\right) \right\Vert _{%
\mathcal{F}}\leq M\left\Vert x\right\Vert _{\mathcal{F}},\forall x\in
X\right\} .
\end{equation*}%
Furthermore, in both cases, denote $\left\Vert \left\vert A\right\vert
\right\Vert _{\mathcal{F}}=\underset{M}{\inf }\mathcal{M}_{A}.$

We have the following.

\textbf{Theorem 3.3} \textit{$\left( i\right) $ If $A:X\mathbb{\rightarrow }%
\mathbf{R}$ is linear and continuous at $\widetilde{0},$ then 
\begin{equation*}
\left| A\left( x\right) \right| \leq \left| \left\| A\right\| \right| _{%
\mathcal{F}}\left\| x\right\| _{\mathcal{F}}
\end{equation*}
for all $x\in X$ and 
\begin{equation*}
\left| \left\| A\right\| \right| _{\mathcal{F}}=\sup \left\{ \left| A\left(
x\right) \right| ;x\in X\text{,}\left\| x\right\| _{\mathcal{F}}\leq
1\right\} .
\end{equation*}
}

\textit{$\left( ii\right) $ If $A:$ $X\mathbb{\rightarrow }$ $X$ is linear
and continuous at $\widetilde{0},$ then 
\begin{equation}
\left\| A\left( x\right) \right\| _{\mathcal{F}}\leq \left| \left\|
A\right\| \right| _{\mathcal{F}}\left\| x\right\| _{\mathcal{F}}  \tag{$1$}
\end{equation}
for all $x\in X$ and 
\begin{equation}
\left| \left\| A\right\| \right| _{\mathcal{F}}=\sup \left\{ \left\| A\left(
x\right) \right\| _{\mathcal{F}};x\in X\text{,}\left\| x\right\| _{\mathcal{F%
}}\leq 1\right\} .  \tag{$2$}
\end{equation}
}

\textbf{Proof.} The proof follows by standard techniques from functional
analysis. For example, let us prove $\left( ii\right) $. Let us suppose that 
$\left( 1\right) $ does not hold, i.e., $\exists x^{^{\prime }}\in X$ with $%
\left\Vert A\left( x^{^{\prime }}\right) \right\Vert _{\mathcal{F}%
}>\left\vert \left\Vert A\right\Vert \right\vert _{\mathcal{F}}\left\Vert
x^{^{\prime }}\right\Vert _{\mathcal{F}},$ which means 
\begin{equation*}
D\left( \widetilde{0},A\left( x^{^{\prime }}\right) \right) >\left\vert
\left\Vert A\right\Vert \right\vert _{\mathcal{F}}D\left( \widetilde{0}%
,x^{^{\prime }}\right) .
\end{equation*}%
Denote $\alpha =D\left( \widetilde{0},x^{^{\prime }}\right) .$ We have $%
\alpha >0$ (otherwise $x^{^{\prime }}=\widetilde{0},$ i.e. $A\left(
x^{^{\prime }}\right) =\widetilde{0}$ and $\left\Vert A\left( \widetilde{0}%
\right) \right\Vert _{\mathcal{F}}=0>\left\vert \left\Vert A\right\Vert
\right\vert _{\mathcal{F}}\left\Vert \widetilde{0}\right\Vert _{\mathcal{F}%
}, $ which is impossible). We get $D\left( \widetilde{0},\frac{1}{\alpha }%
\odot A\left( x^{^{\prime }}\right) \right) >\left\vert \left\Vert
A\right\Vert \right\vert _{\mathcal{F}},$ i.e. $D\left( \widetilde{0}%
,A\left( x^{^{\prime \prime }}\right) \right) >\left\vert \left\Vert
A\right\Vert \right\vert _{\mathcal{F}}$ with $x^{^{\prime \prime }}=\frac{1%
}{\alpha }\odot x^{^{\prime }},\left\Vert x^{^{\prime \prime }}\right\Vert _{%
\mathcal{F}}=1.$ For $\varepsilon =\left\Vert A\left( x^{^{\prime \prime
}}\right) \right\Vert _{\mathcal{F}}-\left\vert \left\Vert A\right\Vert
\right\vert _{\mathcal{F}}>0, $ by $\left\vert \left\Vert A\right\Vert
\right\vert _{\mathcal{F}}=\inf \mathcal{M}_{A},$ there exists $M\in 
\mathcal{M}_{A}$ with $M<\left\vert \left\Vert A\right\Vert \right\vert _{%
\mathcal{F}}+\varepsilon =\left\Vert A\left( x^{^{\prime \prime }}\right)
\right\Vert _{\mathcal{F}},$ which implies the contradiction 
\begin{equation*}
\left\Vert A\left( x^{^{\prime \prime }}\right) \right\Vert _{\mathcal{F}%
}\leq M\left\Vert x^{^{\prime \prime }}\right\Vert _{\mathcal{F}%
}=M<\left\Vert A\left( x^{^{\prime \prime }}\right) \right\Vert _{\mathcal{F}%
}..
\end{equation*}%
Therefore, $\left( 1\right) $ must hold.

It remains to prove $\left( 2\right) .$ First, for $x\in X$, $\left\Vert
x\right\Vert _{\mathcal{F}}\leq 1,$ by $\left( 1\right) $ we get $\left\Vert
A\left( x\right) \right\Vert _{\mathcal{F}}\leq \left\vert \left\Vert
A\right\Vert \right\vert _{\mathcal{F}},$ which implies 
\begin{equation*}
\sup \left\{ \left\Vert A\left( x\right) \right\Vert _{\mathcal{F}};x\in X%
\text{,}\left\Vert x\right\Vert _{\mathcal{F}}\leq 1\right\} \leq \left\vert
\left\Vert A\right\Vert \right\vert _{\mathcal{F}}.
\end{equation*}%
If $\left\vert \left\Vert A\right\Vert \right\vert _{\mathcal{F}}=0,$ this
inequality becomes equality. So, let us assume $\left\vert \left\Vert
A\right\Vert \right\vert _{\mathcal{F}}>0.$ There exists $n_{0}\in \mathbb{N}
$ such that $\left\vert \left\Vert A\right\Vert \right\vert _{\mathcal{F}}-%
\frac{1}{n}>0,$ for all $n\geq n_{0}.$ Since $\left\vert \left\Vert
A\right\Vert \right\vert _{\mathcal{F}}-\frac{1}{n}\notin \mathcal{M}_{A},$
there exists $x^{^{\prime }}\in X$, $x^{^{\prime }}\neq \widetilde{0}$ with 
\begin{equation*}
\left\Vert A\left( x^{^{\prime }}\right) \right\Vert _{\mathcal{F}}>\left(
\left\vert \left\Vert A\right\Vert \right\vert _{\mathcal{F}}-\frac{1}{n}%
\right) \left\Vert x^{^{\prime }}\right\Vert _{\mathcal{F}},
\end{equation*}%
which can be written as $\left\Vert A\left( x^{^{\prime \prime }}\right)
\right\Vert _{\mathcal{F}}>\left( \left\vert \left\Vert A\right\Vert
\right\vert _{\mathcal{F}}-\frac{1}{n}\right) ,$ for all $n\geq n_{0},$ with 
\begin{equation*}
x^{^{\prime \prime }}=\frac{1}{\left\Vert x^{^{\prime }}\right\Vert _{%
\mathcal{F}}}\odot x^{^{\prime }}.
\end{equation*}%
Therefore, 
\begin{equation*}
\sup \left\{ \left\Vert A\left( x\right) \right\Vert _{\mathcal{F}};x\in X%
\text{,}\left\Vert x\right\Vert _{\mathcal{F}}\leq 1\right\} \geq \left\vert
\left\Vert A\right\Vert \right\vert _{\mathcal{F}}-\frac{1}{n}.
\end{equation*}%
Passing to the limit as $n\rightarrow \infty ,$ we get $\sup \left\{
\left\Vert A\left( x\right) \right\Vert _{\mathcal{F}};x\in X\text{,}%
\left\Vert x\right\Vert _{\mathcal{F}}\leq 1\right\} \geq \left\vert
\left\Vert A\right\Vert \right\vert _{\mathcal{F}}$, which proves $\left(
2\right) .\Box $

\textbf{Corollary 3.4} \textit{$\left( i\right) $ If $A:$ $X\mathbb{%
\rightarrow }$ $\mathbf{R}$ is additive (i.e., $A\left( x\oplus y\right)
=A\left( x\right) +A\left( y\right) $), positive homogeneous (i.e., $A\left(
\lambda \odot y\right) =\lambda A\left( x\right) ,\forall \lambda \geq 0$)
and continuous at $\widetilde{0},$ then 
\begin{equation*}
\left\vert A\left( x\right) \right\vert \leq \left\vert \left\Vert
A\right\Vert \right\vert _{\mathcal{F}}\left\Vert x\right\Vert _{\mathcal{F}%
},\forall x\in X\text{.}
\end{equation*}%
$\left( ii\right) $ If $A:$ $X\mathbb{\rightarrow }$ $X$ is additive,
positive homogeneous and continuous at $\widetilde{0},$ then 
\begin{equation*}
||A\left( x\right) ||_{\mathcal{F}}\leq \left\vert \left\Vert A\right\Vert
\right\vert _{\mathcal{F}}\left\Vert x\right\Vert _{\mathcal{F}},\forall
x\in X\text{.}
\end{equation*}%
}

The proof is similar to that of Theorem $3.3,$ because it uses only the
positive homogeneity and additivity (throughout the proof of Theorem $3.2$)
of $A.$

\textbf{Remark.} Examples of operators are the following. Let $X\mathbb{=}%
\mathbf{R}_{\mathcal{F}}.$ Define the following operators $A_{1},A_{4},A_{5}:%
\mathbf{R}_{\mathcal{F}}\rightarrow \mathbf{R},$ $A_{2},A_{3}:\mathbf{R}_{%
\mathcal{F}}\rightarrow \mathbf{R}_{\mathcal{F}}$ by the following
expressions:%
\begin{align*}
A_{1}\left( x\right) & =\int\limits_{0}^{1}\left[ x_{-}\left( r\right)
+x_{+}\left( r\right) \right] dr, \\
A_{2}\left( x\right) & =\left\{ \int\limits_{0}^{1}\left[ x_{+}\left(
0\right) -x_{+}\left( r\right) \right] dr\right\} \odot c, \\
A_{3}\left( x\right) & =\left\{ \int\limits_{0}^{1}\left[ x_{-}\left(
1\right) -x_{-}\left( r\right) \right] dr\right\} \odot c,c\in \mathbf{R}_{%
\mathcal{F}}, \\
A_{4}\left( x\right) & =\int\limits_{0}^{1}x_{-}\left( r\right) dr, \\
A_{5}\left( x\right) & =\int\limits_{0}^{1}x_{+}\left( r\right) dr,
\end{align*}%
where $\left[ x_{-}\left( r\right) ,x_{+}\left( r\right) \right] =\left\{
t\in \mathbf{R};x\left( t\right) \geq r\right\} $. By e.g. $\left[ 11\right]
,$ $x_{-}\left( r\right) $ is bounded nondecreasing on $\left[ 0,1\right] ,$ 
$x_{+}\left( r\right) $ is bounded non-increasing on $\left[ 0,1\right] ,$
both are left continuous on $\left( 0,1\right] $ and right continuous at $%
r=0.$ It follows that $x_{+}\left( 0\right) -x_{+}\left( r\right) \geq
0,x_{-}\left( 1\right) -x_{-}\left( r\right) \geq 0$ and by a simple
calculation, $A_{1}$ is linear and continuous at each $x\in \mathbf{R}_{%
\mathcal{F}},$ and $A_{4},$ $A_{5}$ are additive, positive homogeneous and
continuous at each $x\in \mathbf{R}_{\mathcal{F}}.$ Finally, $A_{2}$ and $%
A_{3}$ are linear continuous operators on $\mathbf{R}_{\mathcal{F}}$. Other
examples of linear operators induced by some fuzzy differential equations
will be considered in the next section.

Next, let us denote 
\begin{align*}
\mathcal{L}_{0}^{+}\left( X\right) & =\left\{ A:X\mathbb{\rightarrow }X;A%
\text{ is additive, positive homogeneous and continuous at }\widetilde{0}%
\right\} , \\
\mathcal{L}_{0}\left( X\right) & =\left\{ A:X\mathbb{\rightarrow }X;A\text{
is linear and continuous at }\widetilde{0}\right\} ,
\end{align*}%
where $\left( X,d\right) $ is any of the spaces described in the beginning
of this section. Let us consider the metric $\Phi :\mathcal{L}_{0}^{+}\left(
X\right) \times \mathcal{L}_{0}^{+}\left( X\right) \rightarrow \mathbf{R}_{%
\mathcal{+}}$ by 
\begin{equation*}
\Phi \left( A,B\right) =\sup \left\{ d\left( A\left( x\right) ,B\left(
x\right) \right) ;\left\Vert x\right\Vert _{\mathcal{F}}\leq 1\right\} .
\end{equation*}%
Clearly, we have $\Phi \left( A,\widetilde{O}\right) =\left\Vert \left\vert
A\right\vert \right\Vert _{\mathcal{F}},A\in \mathcal{L}_{0}^{+}\left(
X\right) ,$ where $\widetilde{O}:X\rightarrow X$ is given by $\widetilde{O}%
\left( x\right) =\widetilde{0},\forall x\in X$.

\textbf{Theorem 3.5} \textit{$\left( \mathcal{L}_{0}^{+}\left( X\right)
,\Phi \right) $ is a complete metric space and, in addition, $\Phi $ has the
following properties: if we define}%
\begin{equation*}
\mathit{\left( A\oplus B\right) \left( x\right) =A\left( x\right) \oplus
B\left( x\right) }
\end{equation*}%
\textit{and }%
\begin{equation*}
\mathit{\left( \lambda \odot A\right) \left( x\right) =\lambda \odot A\left(
x\right) ,}
\end{equation*}%
\textit{then the following hold:}

(i) $\Phi \left( A\oplus B,C\oplus D\right) \leq \Phi \left( A,C\right)
+\Phi \left( B,D\right) ,$

(ii) $\Phi \left( k\odot A,k\odot B\right) =\left\vert k\right\vert \Phi
\left( A,B\right) ,$

(iii) $\Phi \left( A,B\right) \leq \left\vert \left\Vert A\right\Vert
\right\vert _{\mathcal{F}}+\left\vert \left\Vert B\right\Vert \right\vert _{%
\mathcal{F}},$

(iv) $\Phi \left( A\oplus B,C\right) \leq \Phi \left( A,C\right) +\Phi
\left( B,C\right) ,$

(v) $\Phi \left( A\oplus B,\widetilde{O}\right) \leq \left\vert \left\Vert
A\right\Vert \right\vert _{\mathcal{F}}+\left\vert \left\Vert B\right\Vert
\right\vert _{\mathcal{F}}.$

\textbf{Proof.} The properties of $\Phi $ are immediate consequences of the
properties of $d.$ Let us prove that $\left( \mathcal{L}_{0}^{+}\left(
X\right) ,\Phi \right) $ is complete. Let $A_{n}\in \mathcal{L}%
_{0}^{+}\left( X\right) ,n\in \mathbb{N}$ be a Cauchy sequence, i.e. $%
\forall \varepsilon >0,\exists n_{0}\in \mathbb{N}$ such that $\Phi \left(
A_{n},A_{n+p}\right) <\varepsilon ,$ for all $n\geq n_{0},p\in \mathbb{N}.$
For $x\in X$, we will define $A:X\mathbb{\rightarrow }X$ as follows. If $%
\left\Vert x\right\Vert _{\mathcal{F}}\leq 1,$ then by the definition of $%
\Phi $ we get 
\begin{equation*}
d\left( A_{n}\left( x\right) ,A_{n+p}\left( x\right) \right) <\varepsilon ,%
\text{for all }n\geq n_{0},p\in \mathbb{N}\text{.}
\end{equation*}%
Moreover, if $\left\Vert x\right\Vert _{\mathcal{F}}>1,$ denoting $y=\alpha
\odot x$ with $\alpha =\frac{1}{\left\Vert x\right\Vert _{\mathcal{F}}}>0,$
we get $\left\Vert y\right\Vert _{\mathcal{F}}=1,d\left( A_{n}\left(
y\right) ,A_{n+p}\left( y\right) \right) =\alpha d\left( A_{n}\left(
x\right) ,A_{n+p}\left( x\right) \right) <\varepsilon ,$ i.e., 
\begin{equation*}
d\left( A_{n}\left( x\right) ,A_{n+p}\left( x\right) \right) <\varepsilon
\left\Vert x\right\Vert _{\mathcal{F}},\text{ for all }n\geq n_{0},p\in 
\mathbb{N}.
\end{equation*}%
Consequently, $\left( A_{n}\left( x\right) \right) _{n\in \mathbb{N}}$ is a
Cauchy sequence in the complete metric space $\left( X,d\right) ,$ i.e., it
is convergent. Let us denote $A\left( x\right) =\underset{n\rightarrow
\infty }{\lim }A_{n}\left( x\right) ,$ i.e.,%
\begin{equation*}
\underset{n\rightarrow \infty }{\lim }d\left( A_{n}\left( x\right) ,A\left(
x\right) \right) =0,x\in X.
\end{equation*}%
We have to prove that $A\in \mathcal{L}_{0}^{+}\left( X\right) .$ First, 
\begin{align*}
0& \leq d\left( A\left( x\oplus y\right) ,A\left( x\right) \oplus A\left(
y\right) \right) \leq d\left( A\left( x\oplus y\right) ,A_{n}\left( x\oplus
y\right) \right) \\
& +d\left( A_{n}\left( x\oplus y\right) ,A_{n}\left( x\right) \oplus
A_{n}\left( y\right) \right) +d\left( A_{n}\left( x\right) \oplus
A_{n}\left( y\right) ,A\left( x\right) \oplus A\left( y\right) \right) \\
& \leq d\left( A\left( x\oplus y\right) ,A_{n}\left( x\oplus y\right)
\right) +d\left( A_{n}\left( x\right) ,A\left( x\right) \right) +d\left(
A_{n}\left( y\right) ,A\left( y\right) \right) .
\end{align*}%
Passing to the limit as $n\rightarrow \infty ,$ we obtain $0\leq d\left(
A\left( x\oplus y\right) ,A\left( x\right) \oplus A\left( y\right) \right)
\leq 0,$ i.e., $A\left( x\oplus y\right) =A\left( x\right) \oplus A\left(
y\right) .$ Let $\lambda >0$ and $x\in X$. Similarly, we get 
\begin{align*}
0& \leq d\left( A\left( \lambda \odot x\right) ,\lambda \odot A\left(
x\right) \right) \leq d\left( A\left( \lambda \odot x\right) ,A_{n}\left(
\lambda \odot x\right) \right) \\
& +d\left( A_{n}\left( \lambda \odot x\right) ,\lambda \odot A_{n}\left(
x\right) \right) +d\left( \lambda \odot A_{n}\left( x\right) ,\lambda \odot
A\left( x\right) \right) \\
& \leq d\left( A\left( \lambda \odot x\right) ,A_{n}\left( \lambda \odot
x\right) \right) +\lambda d\left( A_{n}\left( x\right) ,A\left( x\right)
\right) ,
\end{align*}%
and passing to the limit as $n\rightarrow \infty ,$ it follows that $A\left(
\lambda \odot x\right) =\lambda \odot A\left( x\right) .$

Now, we show that $A$ is continuous at $\widetilde{0}.$ First, by Corollary $%
3.4,\left( ii\right) ,$ we have $\left\Vert A_{n}\left( x\right) \right\Vert
_{\mathcal{F}}\leq \left\vert \left\Vert A_{n}\right\Vert \right\vert _{%
\mathcal{F}}\left\Vert x\right\Vert _{\mathcal{F}},\forall x\in X$, $n\in 
\mathbb{N}$. We have 
\begin{align*}
\left\Vert A\left( x\right) \right\Vert _{\mathcal{F}}& =d\left( A\left(
x\right) ,\widetilde{0}\right) \leq d\left( A\left( x\right) ,A_{n}\left(
x\right) \right) +d\left( A_{n}\left( x\right) ,\widetilde{0}\right) \\
& \leq d\left( A\left( x\right) ,A_{n}\left( x\right) \right) +\left\vert
\left\Vert A_{n}\left( x\right) \right\Vert \right\vert _{\mathcal{F}%
}\left\Vert x\right\Vert _{\mathcal{F}}.
\end{align*}%
But, since $\Phi $ is a metric, we obtain 
\begin{equation*}
\left\vert \Phi \left( A_{n},\widetilde{O}\right) -\Phi \left( A_{m},%
\widetilde{O}\right) \right\vert \leq \Phi \left( A_{n},A_{m}\right) ,
\end{equation*}%
which shows that the sequence of real positive numbers $\left( \Phi \left(
A_{n},\widetilde{O}\right) \right) _{n\in \mathbb{N}}$ is a Cauchy sequence,
i.e., it is bounded, implying the existence of $M>0$ with 
\begin{equation*}
\left\vert \left\Vert A_{n}\left( x\right) \right\Vert \right\vert _{%
\mathcal{F}}=\Phi \left( A_{n},\widetilde{O}\right) \leq M,\forall n\in 
\mathbb{N}\text{.}
\end{equation*}%
It follows that 
\begin{equation*}
\left\Vert A\left( x\right) \right\Vert _{\mathcal{F}}\leq d\left( A\left(
x\right) ,A_{n}\left( x\right) \right) +M\left\Vert x\right\Vert _{\mathcal{F%
}},\forall x\in X\text{, }n\in \mathbb{N}\text{.}
\end{equation*}%
Passing to the limit as $n\rightarrow \infty ,$ we obtain 
\begin{equation*}
\left\Vert A\left( x\right) \right\Vert _{\mathcal{F}}\leq M\left\Vert
x\right\Vert _{\mathcal{F}},\forall x\in X\text{,}
\end{equation*}%
i.e., according to Theorem $3.2$ $\left( ii\right) ,$ $A$ is continuous at $%
\widetilde{0}\in X$. In conclusion, $\left( \mathcal{L}_{0}^{+}\left(
X\right) ,\Phi \right) $ is complete, which proves the theorem. $\Box $

\textbf{Corollary 3.6} \textit{$\left( i\right) $ $\left( \mathcal{L}%
_{0}\left( X\right) ,\Phi\right) $ is a complete metric space. }

\textit{$\left( ii\right) $ If we denote 
\begin{align*}
\mathcal{L}^{+}\left( X\right) & =\left\{ A\in \mathcal{L}_{0}^{+}\left(
X\right) ;A\text{ is continuous at each }x\in X\right\} , \\
\mathcal{L}\left( X\right) & =\left\{ A\in \mathcal{L}_{0}\left( X\right) ;A%
\text{ is continuous at each }x\in X\right\} ,
\end{align*}%
then $\left( \mathcal{L}^{+}\left( X\right) ,\Phi \right) $ and $\left( 
\mathcal{L}\left( X\right) ,\Phi \right) $ are complete metric spaces.}

Concerning these spaces of operators, in what follows, we prove the uniform
boundedness principle.

\textbf{Theorem 3.7} \textit{Let $\left( X,d\right) $ be any of the spaces
listed in the beginning of this section, and let $\mathbf{L}\left( X\right) $
be either $\mathcal{L}\left( X\right) $ or $\mathcal{L}^{+}\left( X\right) $%
. If $A_{j}\in \mathbf{L}\left( X\right) ,j\in J,$ is pointwise bounded,
i.e., for any $x\in X$, $\left\Vert A_{j}\left( x\right) \right\Vert _{%
\mathcal{F}}=d\left( A_{j}\left( x\right) ,\widetilde{0}\right) \leq M_{x},$
for all $j\in J,$ then there exists $M>0$ such that 
\begin{equation*}
\left\vert \left\Vert A_{j}\right\Vert \right\vert _{\mathcal{F}}\leq
M,\forall j\in J\text{,}
\end{equation*}%
(i.e., $\left( A_{j}\right) _{j}$ is uniformly bounded).}

\textbf{Proof.} For any $n\in \mathbb{N}$, let us denote $A_{n}=\left\{ x\in
X;\text{ }\left\| A_{j}\left( x\right) \right\| _{\mathcal{F}}\leq n,\forall
j\in J\right\} .$ It is obvious that $X\mathbb{=}\underset{n\in \mathbb{N}}{%
\bigcup }A_{n}.$ But $A_{n}$ are closed sets (because if $d\left(
x_{m},x\right) \overset{m\rightarrow \infty }{\longrightarrow }0,x_{m}\in
A_{n},\forall m,$ then by $d\left( A_{j}\left( x\right) ,\widetilde{0}%
\right) \leq d\left( A_{j}\left( x\right) ,A_{j}\left( x_{m}\right) \right)
+d\left( A_{j}\left( x_{m}\right) ,\widetilde{0}\right) \leq n+d\left(
A_{j}\left( x\right) ,A_{j}\left( x_{m}\right) \right) ,$ passing to the
limit as $m\rightarrow \infty $ and taking into account the continuity of $%
A_{j}$ at each $x,$ we get $x\in A_{n}$).

Since $\left( X,d\right) $ is a complete metric space, it is of second Baire
category, therefore there exists $m\in \mathbb{N}$, such that $intA_{m}\neq
\emptyset .$ Let $x_{0}\in intA_{m}$ and $\lambda >0$ such that $B\left(
x_{0},\lambda \right) =\left\{ x\in X\text{; }d\left( x,\widetilde{0}\right)
<\lambda \right\} \subset A_{m}.$ For $x\in X$, denote $x_{1}=x_{0}\oplus 
\frac{\lambda }{2\left\Vert x\right\Vert _{\mathcal{F}}}\odot x.$ We have 
\begin{align*}
d\left( x_{1},x_{0}\right) & =d\left( x_{0}\oplus \frac{\lambda }{%
2\left\Vert x\right\Vert _{\mathcal{F}}}\odot x,x_{0}\right) =d\left( \frac{%
\lambda }{2\left\Vert x\right\Vert _{\mathcal{F}}}\odot x,\widetilde{0}%
\right) \\
& =\frac{\lambda }{2\left\Vert x\right\Vert _{\mathcal{F}}}d\left( x,%
\widetilde{0}\right) =\frac{\lambda }{2}<\lambda ,
\end{align*}%
i.e., $x_{1}\in B\left( x_{0},\lambda \right) .$ Then, for all $x\in X$, $%
x\neq \widetilde{0},$ we have 
\begin{align*}
d\left( A_{j}\left( x\right) ,\widetilde{0}\right) & =\frac{2\left\Vert
x\right\Vert _{\mathcal{F}}}{\lambda }d\left( A_{j}\left( x_{0}\oplus \frac{%
\lambda }{2\left\Vert x\right\Vert _{\mathcal{F}}}\odot x\right)
,A_{j}\left( x_{0}\right) \right) \\
& =\frac{2\left\Vert x\right\Vert _{\mathcal{F}}}{\lambda }d\left(
A_{j}\left( x_{1}\right) ,A_{j}\left( x_{0}\right) \right) \leq \frac{%
2\left\Vert x\right\Vert _{\mathcal{F}}}{\lambda }\left[ \left\Vert
A_{j}\left( x_{1}\right) \right\Vert _{\mathcal{F}}+\left\Vert A_{j}\left(
x_{0}\right) \right\Vert _{\mathcal{F}}\right] \\
& \leq \frac{4m}{\lambda }\left\Vert x\right\Vert _{\mathcal{F}},
\end{align*}%
which implies $\left\Vert A_{j}\left( x\right) \right\Vert _{\mathcal{F}%
}\leq \frac{4m}{\lambda }\left\Vert x\right\Vert _{\mathcal{F}}$ and,
therefore, $\left\vert \left\Vert A_{j}\left( x\right) \right\Vert
\right\vert _{\mathcal{F}}\leq \frac{4m}{\lambda }=M,$ for all $j\in J.$ The
theorem is proved. $\Box $

\textbf{Remark.} Some results in classical functional analysis concerning
invertible operators can also be considered but only in a particular case.

Thus, let us recall that $x\in \mathbf{R}_{\mathcal{F}}$ is called
triangular fuzzy number, if we have $\left[ x\right] ^{1}=\left\{
x_{c}\right\} ,\left[ x\right] ^{0}=\left[ x_{l},x_{r}\right] $ where $%
x_{l}\leq x_{c}\leq x_{r}$ (see e.g. $\left[ 8\right] $). This implies%
\begin{equation*}
\left[ x\right] ^{r}=\left[ x_{c}-\left( 1-r\right) \left(
x_{c}-x_{r}\right) ,x_{c}+\left( 1-r\right) \left( x_{r}-x_{c}\right) \right]
,\text{ }r\in \left[ 0,1\right] .
\end{equation*}%
We denote $x=\left( x_{l},x_{c},x_{r}\right) .$ A triangular number $x$ is
called symmetric if there exists $\delta \geq 0$ with $x=\left( x_{c}-\delta
,x_{c},x_{c}+\delta \right) .$ Denote by $\mathbf{R}_{\mathcal{F}}^{TS}$ the
set of all symmetric triangular fuzzy numbers. Under $\oplus ,\odot $ and
passing to limit, $\mathbf{R}_{\mathcal{F}}^{TS}$ is closed, therefore $%
\left( \mathbf{R}_{\mathcal{F}}^{TS},D\right) $ is a complete metric space.
Also, for any $x_{1},x_{2}\in $ $\mathbf{R}_{\mathcal{F}}^{TS},$ there
exists $x_{1}\ominus x_{2}\in $ $\mathbf{R}_{\mathcal{F}}^{TS}$ or $%
x_{2}\ominus x_{1}\in $ $\mathbf{R}_{\mathcal{F}}^{TS}$ (see e.g. $\left[ 7%
\right] ,\left[ 17\right] $).

This last property allows us to state the following

\textbf{Theorem 3.8} \textit{Let $A\in \mathcal{L}\left( \mathbf{R}_{%
\mathcal{F}}^{TS}\right) .$ Then $A^{-1}\in \mathcal{L}\left( \mathbf{R}_{%
\mathcal{F}}^{TS}\right) $ if and only if there exists $\lambda >0$ such
that 
\begin{equation*}
\left\| A\left( x\right) \right\| _{\mathcal{F}}\geq \lambda \left\|
x\right\| _{\mathcal{F}},\forall x\in \mathbf{R}_{\mathcal{F}}^{TS}.
\end{equation*}
In this case, $\left\| \left| A^{-1}\right| \right\| _{\mathcal{F}}\leq 
\frac{1}{\lambda }.$}

\textbf{Proof.} First, let us assume $A^{-1}\in \mathcal{L}\left( \mathbf{R}%
_{\mathcal{F}}^{TS}\right) .$ Similarly (see the proof of Theorem $%
3.2.,\left( ii\right) $), there exists $\alpha >0$ such that 
\begin{equation*}
\left\Vert A^{-1}\left( y\right) \right\Vert _{\mathcal{F}}\leq \alpha
\left\Vert y\right\Vert _{\mathcal{F}},\forall y\in A\left( \mathbf{R}_{%
\mathcal{F}}^{TS}\right) .
\end{equation*}%
For $x=A^{-1}\left( y\right) ,$ we get $\left\Vert A\left( x\right)
\right\Vert _{\mathcal{F}}\geq \frac{1}{\alpha }\left\Vert x\right\Vert _{%
\mathcal{F}},$ i.e., we can take $\lambda =\frac{1}{\alpha }.$

Conversely, let us suppose that $\left\Vert A\left( x\right) \right\Vert _{%
\mathcal{F}}\geq \lambda \left\Vert x\right\Vert _{\mathcal{F}},\forall x\in 
\mathbf{R}_{\mathcal{F}}^{TS}$ and take $x_{1},x_{2}\in \mathbf{R}_{\mathcal{%
F}}^{TS},x_{1}\neq x_{2}.$ Then, $x_{1}\ominus x_{2}$ or $x_{2}\ominus x_{1}$
exists, so let us choose, for example, $u=x_{1}\ominus x_{2},$ i.e., $%
x_{1}=x_{2}\oplus u.$ We get $A\left( x_{1}\right) \ominus A\left(
x_{2}\right) =A\left( u\right) .$ Now, if $A\left( x_{1}\right) =A\left(
x_{2}\right) ,$ we obtain $A\left( u\right) =\widetilde{0}$ and by our
inequality, it follows $u=\widetilde{0},$ i.e., $x_{1}=x_{2}.$ Therefore $A$
is invertible and, moreover, for $y=A\left( x\right) ,$%
\begin{equation*}
\left\Vert A^{-1}\left( y\right) \right\Vert _{\mathcal{F}}=\left\Vert
x\right\Vert _{\mathcal{F}}\leq \frac{1}{\lambda }\left\Vert A\left(
x\right) \right\Vert _{\mathcal{F}}=\frac{1}{\lambda }\left\Vert
y\right\Vert _{\mathcal{F}},
\end{equation*}%
which implies $A^{-1}\in \mathcal{L}\left( \mathbf{R}_{\mathcal{F}%
}^{TS}\right) $ and $\left\Vert \left\vert A^{-1}\right\vert \right\Vert _{%
\mathcal{F}}\leq \frac{1}{\lambda }.$ $\Box $

\textbf{Remark.} Unfortunately, if we replace $\mathbf{R}_{\mathcal{F}}^{TS}$
by any from the function spaces $X=\left( L^{p}\left( \left[ a,b\right] ;%
\mathbf{R}_{\mathcal{F}}\right) ,D_{p}\right) ,\left( C^{p}\left( \left[ a,b%
\right] ;\mathbf{R}_{\mathcal{F}}\right) ,D_{p}^{\ast }\right) ,p\in \mathbb{%
N}$, Theorem $3.8$ does not hold because we cannot define, in general, $%
x_{1}\ominus x_{2}$ and $x_{2}\ominus x_{1}$ for $x_{1},x_{2}\in X$.
Effectivelly, this means that the uniqueness problem in fuzzy differential
equations may not be valid in these spaces. At the end of this section, we
consider some elements of semigroup theory.

\textbf{Theorem 3.9} \textit{Let $\left( X,d\right) $ be any from the spaces
described in the beginning of this section, and $\left( L\left( X\right)
,\Phi \right) $ be the space $\left( \mathcal{L}\left( X\right) ,\Phi
\right) .$ Let us define $T\left( t\right) =e^{t\odot A},$ $t\in $}$\mathbf{R%
}$,\textit{\ by%
\begin{equation*}
\underset{m\rightarrow \infty }{\lim }\Phi \left( T\left( t\right)
,\sum\limits_{p=0}^{m}\frac{t^{p}}{p!}\odot A^{p}\right) =0,
\end{equation*}%
where $\sum $ is the sum with respect to $\oplus ,A\in L\left( X\right) $
and $A^{0}=I,A^{p}=A^{p-1}\circ A,p=2,3...,.$ Formally, we write 
\begin{equation*}
e^{t\odot A}=\sum\limits_{p=0}^{\infty }\frac{t^{p}}{p!}\odot A^{p}.
\end{equation*}%
We have:}

\textit{$\left( i\right) $ $T\left( t\right) \in L\left( X\right) $ for all $%
t\in \mathbf{R};$ }

\textit{$\left( ii\right) $ $T\left( t+s\right) =T\left( t\right)(T\left(
s\right)), $ for all $t,s\geq 0, $ or for all $t,s\le 0$; }

\textit{(iii) There exists }%
\begin{equation*}
\mathit{\lim_{h\searrow 0}\frac{1}{h}\odot \lbrack T(h)(x)\ominus x]=A(x),}
\end{equation*}%
\textit{for all $x\in X$, where the limit is considered in the metric d;}

\textit{$\left( iv\right) $ $T\left( t\right) $ is continuous as function of 
$t\in $}$\mathbf{R}$\textit{\ and $T\left( 0\right) =I.$ Also, $T(t)$ is
generalized differentiable with respect to all $t\in $}$\mathbf{R}$\textit{,
with the derivative equal to $A[T(t)]$. More exactly, it is Hukuhara
differentiable (in the sense of Definition 2.7, (i), see the remark at the
end of Section 2) with respect to $t\in \mathbf{R}_{+}$ , i.e., for all $%
t\geq 0$ we have 
\begin{align*}
& \lim_{h\searrow 0}d(\frac{1}{h}\odot (T(t+h)(x)\ominus T(t)(x)),A[T(t)(x)])
\\
& =\lim_{h\searrow 0}d(\frac{1}{h}\odot (T(t)(x)\ominus
T(t-h)(x)),A[T(t)(x)]) \\
& =0,
\end{align*}%
and generalized differentiable (in the sense of Definition 2.7,(ii)) with
respect to $t<0$, i.e., for all $t<0$ we have 
\begin{align*}
& \lim_{h\searrow 0}d(-\frac{1}{h}\odot (T(t)(x)\ominus
T(t+h)(x)),A[T(t)(x)]) \\
& =\lim_{h\searrow 0}d(-\frac{1}{h}\odot (T(t-h)(x)\ominus
T(t)(x)),A[T(t)(x)]) \\
& =0,
\end{align*}%
for all $x\in X$. }

\textit{(v) If $u_{0}\in X$ and $g:\mathbf{R}\rightarrow X$ is continuous on 
}$\mathbf{R}$\textit{, then 
\begin{equation*}
u(t)=T(t)(u_{0})\oplus \int_{0}^{t}T(t-s)g(s)ds
\end{equation*}%
is generalized differentiable on }$\mathbf{R}$\textit{\ (more exactly, it is
Hukuhara differentiable on $\mathbf{R}_{+}$ and generalized differentiable
as in the above point (iv) for $t<0$ ) and satisfies}%
\begin{align*}
\mathit{u^{\prime }(t)}& \mathit{=A[u(t)]\oplus g(t),}\text{ }\mathit{t\in R,%
} \\
u\left( 0\right) & =\mathit{u_{0},}
\end{align*}%
\textit{where $u^{\prime }(t)$ denotes the generalized derivative. Here, the
integral for functions defined on a compact interval with values in $X$ is
considered in the Riemann (classical) sense.}

\textbf{Proof.} (i) Denote $S_{m}\left( t\right) \left( x\right)
=\sum\limits_{p=0}^{m}\frac{t^{p}}{p!}\odot A^{p}(x).$ By Theorem $3.5$ and
Corollary $3.6,$ it suffices to show that $\left( S_{m}\left( t\right)
\right) _{m}$ is a Cauchy sequence in the complete metric space $\left(
L\left( X\right) ,\Phi \right) ,$ for all $t\in \mathbf{R}$. First, since $%
A\in L\left( X\right) $, it follows $A^{p}\in L\left( X\right) $ for $%
p=2,3....$, and $S_{m}\left( t\right) \in L\left( X\right) ,\forall
m=0,1,....$ Then

\begin{align*}
d\left( S_{n}\left( t\right) (x),S_{n+p}\left( t\right) (x)\right) &
=d\left( \sum\limits_{i=0}^{n}\frac{t^{i}}{i!}\odot
A^{i}(x),\sum\limits_{i=0}^{n+p}\frac{t^{i}}{i!}\odot A^{i}(x)\right) \\
& \leq \sum\limits_{i=n+1}^{n+p}d\left( \widetilde{0},\frac{t^{i}}{i!}\odot
A^{i}(x)\right) =\sum\limits_{i=n+1}^{n+p}\frac{|t|^{i}}{i!}d\left( 
\widetilde{0},A^{i}\left( x\right) \right) .
\end{align*}%
However,%
\begin{equation*}
d\left( \widetilde{0},A\left( A\left( x\right) \right) \right) =\left\Vert
A^{2}\left( x\right) \right\Vert _{\mathcal{F}}\leq \left\Vert \left\vert
A\right\vert \right\Vert _{\mathcal{F}}\left\Vert A\left( x\right)
\right\Vert _{\mathcal{F}}\leq \left\Vert \left\vert A\right\vert
\right\Vert _{\mathcal{F}}^{2}\left\Vert x\right\Vert _{\mathcal{F}}
\end{equation*}%
and, so by induction, one can obtain 
\begin{equation*}
\left\Vert \left\vert A^{i}\right\vert \right\Vert _{\mathcal{F}}\leq
\left\Vert \left\vert A\right\vert \right\Vert _{\mathcal{F}}^{i},\text{ for
all }i=2,3,...\mathbf{.}
\end{equation*}%
We obtain%
\begin{equation*}
d\left( S_{n}\left( t\right) (x),S_{n+p}\left( t\right) (x)\right) \leq
\sum\limits_{i=n+1}^{n+p}\frac{\left( |t||||A|||_{\mathcal{F}}\left\Vert
x|\right\Vert _{\mathcal{F}}\right) ^{i}}{i!},
\end{equation*}%
and passing to supremum with $\left\Vert x\right\Vert _{\mathcal{F}}\leq 1,$
we get%
\begin{equation*}
\Phi \left( S_{n}\left( t\right) ,S_{n+p}\left( t\right) \right) \leq
\sum\limits_{i=n+1}^{n+p}\frac{\left( |t|\left\Vert \left\vert A\right\vert
\right\Vert _{\mathcal{F}}\right) ^{i}}{i!},
\end{equation*}%
which immediately implies that $\left( S_{n}\left( t\right) \right) _{n}$ is
a Cauchy sequence, therefore, its limit $T\left( t\right) $ exists in $%
L\left( X\right) $. Note that the last inequality actually implies that $%
\left( S_{n}\left( t\right) \right) _{n}$ is a uniformly Cauchy sequence on $%
[-a,+a]$, for any real number $a>0$, which implies that $\lim_{n\rightarrow
+\infty }\Phi (S_{n}(t),T(t))=0$ holds uniformly on each compact interval $%
[-a,a]$, $a>0$.

(ii) First let $t,s\geq 0$. With the notation in the above point (i), a
simple calculation shows $S_{m}(t)[S_{n}(s)(x)]=S_{m+n}(t+s)(x)$, for all $%
x\in X$. Then 
\begin{align*}
d[T(t)(T(s)(x)),T(t+s)(x)]& \leq d[T(t)(T(s)(x)),S_{m}(t)(T(s)(x))] \\
& +d[S_{m}(t)(T(s)(x)),S_{m}(t)(S_{n}(s)(x))] \\
& +d[S_{m}(t)(S_{n}(s)(x)),S_{m+n}(t+s)(x)] \\
& +d[S_{m+n}(t+s)(x),T(t+s)(x)].
\end{align*}%
Let $\epsilon >0$ be arbitrary fixed. There exists $m_{1},n_{1}\in \mathbf{N}
$, such that for all $m>m_{1}$ and $n>n_{1},$ we have $%
d[S_{m+n}(t+s)(x),T(t+s)(x)]<\epsilon /3$. Because $S_{m}(t)(y)$ converges
to $T(t)(y)$, there exists $m_{2}\in \mathbf{N}$ such that $%
d[T(t)(T(s)(x)),S_{m}(t)(T(s)(x))]<\epsilon /3,\forall m>m_{2}$. Let $%
m>max\{m_{1},m_{2}\}=m_{0}$ be fixed. Since $S_{m}(t)$ is a continuous
linear operator and $S_{n}(s)(x)\rightarrow T(s)(x)$, there exists $n_{2}\in 
\mathbf{N}$ such that for all $n>n_{2}$ we have $%
d[S_{m}(t)(T(s)(x)),S_{m}(t)(S_{n}(s)(x))]<\epsilon /3$. Next, choosing $%
m>m_{0}$ and $n>max\{n_{1},n_{2}\}=n_{0},$ we obtain $%
d[T(t)(T(s)(x)),T(t+s)(x)]<\epsilon /3+\epsilon /3+\epsilon /3=\epsilon ,$
with arbitrary $\epsilon >0$, i.e., $d[T(t)(T(s)(x)),T(t+s)(x)]=0$, which
proves the assertion.

Now, because $t,s\geq 0$, a similar calculation shows that $%
S_{m+n}(-t-s)(x)=S_{m}(-t)(S_{n}(-s)(x))$. Thus, we obtain 
\begin{equation*}
T(t+s)(x)=T(t)(T(s)(x)),\forall t,s\leq 0,
\end{equation*}%
which proves the point (ii).

(iii) Since $S_{m}(h)(x)=x\oplus \sum_{p=1}^{m}\frac{h^{p}}{p!}\odot
A^{p}(x) $, we get 
\begin{equation*}
\frac{1}{h}\odot \lbrack S_{m}(h)(x)\ominus x]=\sum_{p=1}^{m}\frac{h^{p-1}}{%
p!}\odot A^{p}(x).
\end{equation*}%
We need the following auxiliary result in $(X,\oplus ,\odot ,d)$: if $%
d(a_{n},a)\rightarrow 0$ and there exist $a_{n}\ominus b$, for all $n\in 
\mathbf{N}$, then there exists $a\ominus b$ and $d(a_{n}\ominus b,a\ominus
b)\rightarrow 0$. Indeed, denote $c_{n}=a_{n}\ominus b,n\in \mathbf{N}$,
that is $a_{n}=c_{n}\oplus b$, for all $n\in \mathbf{N}$. Since 
\begin{align*}
d(c_{n},c_{m})& =d(a_{n}\ominus b,a_{m}\ominus b) \\
& =d((a_{n}\ominus b)\oplus b,(a_{m}\ominus b)\oplus b) \\
& =d(a_{n},a_{m}),
\end{align*}%
it follows that $(c_{n})_{n}$ is a Cauchy sequence in the complete metric
space $(X,d)$, i.e., it is convergent, so let us denote by $c$ its limit. We
have 
\begin{equation*}
d(a,b\oplus c)\leq d(a,a_{n})+d(c_{n}\oplus b,c\oplus
b)=d(a,a_{n})+d(c_{n},c)\rightarrow 0,
\end{equation*}%
i.e., $a=b\oplus c$, which implies $c=a\ominus b$. Therefore, $%
d(a_{n}\ominus b,a\ominus b)=d(c_{n},c)\rightarrow 0$.

Since $S_{m}(h)(x)\rightarrow T(h)(x)$ when $m\rightarrow +\infty $, for all 
$h>0$ and $x\in X$, applying the auxiliary result, it follows that for all $%
h>0$ and $x\in X,$ we have 
\begin{equation*}
\lim_{m\rightarrow +\infty }d(\frac{1}{h}\odot \lbrack S_{m}(h)(x)\ominus x],%
\frac{1}{h}\odot \lbrack T(h)(x)\ominus x])=0.
\end{equation*}%
Moreover,%
\begin{align*}
& d(\frac{1}{h}\odot \lbrack T(h)(x)\ominus x],A(x)] \\
& \leq d(\frac{1}{h}\odot \lbrack T(h)(x)\ominus x],\frac{1}{h}\odot \lbrack
S_{m}(h)(x)\ominus x]) \\
& +d(\frac{1}{h}\odot \lbrack S_{m}(h)(x)\ominus x],A(x)).
\end{align*}%
Recall that%
\begin{equation*}
\frac{1}{h}\odot \lbrack S_{m}(h)(x)\ominus x]=A(x)\oplus \sum_{p=2}^{m}%
\frac{h^{p-1}}{p!}\odot A^{p}(x),
\end{equation*}%
i.e., there exists $\frac{1}{h}\odot \lbrack S_{m}(h)(x)\ominus x]\ominus
A(x)=\sum_{p=2}^{m}\frac{h^{p-1}}{p!}\odot A^{p}(x)$ and 
\begin{align*}
d(\frac{1}{h}\odot \lbrack S_{m}(h)(x)\ominus x]\ominus A(x),\widetilde{0})&
=d(\frac{1}{h}\odot \lbrack S_{m}(h)(x)\ominus x],A(x)) \\
& =d(\sum_{p=2}^{m}\frac{h^{p-1}}{p!}\odot A^{p}(x),\widetilde{0}) \\
& \leq \sum_{p=2}^{m}\frac{h^{p-1}}{p!}||A^{p}(x)||_{\mathcal{F}} \\
& \leq ||x||_{\mathcal{F}}\sum_{p=2}^{m}\frac{h^{p}}{p!h}|||A|||_{\mathcal{F}%
}^{p} \\
& =||x||_{\mathcal{F}}\frac{e_{m}(h|||A|||_{\mathcal{F}})-1-h|||A|||_{%
\mathcal{F}}}{h},
\end{align*}%
where $e_{m}(h)=1+\frac{h}{1!}+...+\frac{h^{m}}{m!}$ is the partial sum of
order $m$ of the usual exponential $e^{h}$. Consequently, passing to the
limit as $m\rightarrow +\infty $, we obtain 
\begin{equation*}
d(\frac{1}{h}\odot \lbrack T(h)(x)\ominus x],A(x))\leq ||x||_{\mathcal{F}}%
\frac{e^{h|||A|||_{\mathcal{F}}}-1-h|||A|||_{\mathcal{F}}}{h},\forall x\in X.
\end{equation*}%
Now passing to the limit as $h\searrow 0$ and taking into account the
continuity of the metric $d$ (with respect to its components), it follows 
\begin{equation*}
\lim_{h\searrow 0}d(\frac{1}{h}\odot \lbrack T(h)(x)\ominus x],A(x))=0,
\end{equation*}%
for all $x\in X$, which proves (iii).

(iv) We distinguish three cases: $t>0,t=0$ or $t<0$. We have 
\begin{align*}
d\left( T\left( t\right) (x),T\left( t+h\right) (x)\right) & \leq d\left(
T\left( t\right) (x),S_{n}\left( t\right) (x)\right) \\
& +d\left( S_{n}\left( t\right) (x),S_{n}\left( t+h\right) (x)\right) \\
& +d\left( S_{n}\left( t+h\right) (x),T\left( t+h\right) (x)\right) ,
\end{align*}%
where 
\begin{equation*}
S_{n}\left( t\right) (x)=\sum\limits_{p=0}^{n}\frac{t^{p}}{p!}\odot
A^{p}\left( x\right) .
\end{equation*}%
Passing to the limit as $n\rightarrow \infty $ in the above inequality, and
using $\left( viii\right) $ of Theorem $2.4$, we obtain 
\begin{equation*}
d\left( T\left( t\right) (x),T\left( t+h\right) (x)\right) \leq
\sum\limits_{p=1}^{\infty }\frac{\left\vert \left( t+h\right)
^{p}-t^{p}\right\vert }{p!}d\left( \widetilde{0},A^{p}\left( x\right)
\right) .
\end{equation*}%
While for $t=0,$ the above inequality is trivial, when $t>0$ or $t<0$ we
choose a sufficiently small $\delta >0$ such that for all $\left\vert
h\right\vert <\delta ,$ we have either $t\pm h>0$ or $t\pm h<0,$
respectively. Passing to supremum with $\left\Vert x\right\Vert _{\mathcal{F}%
}\leq 1,$ we get 
\begin{equation*}
\Phi \left( T\left( t\right) ,T\left( t+h\right) \right) \leq
\sum\limits_{p=1}^{\infty }\frac{\left\vert \left( t+h\right)
^{p}-t^{p}\right\vert }{p!}\left\Vert \left\vert A\right\vert \right\Vert _{%
\mathcal{F}}^{p}.
\end{equation*}%
But for $t\not=0$, there exists $\sigma $ with $\left\vert \sigma
\right\vert <2|t|,$ such that for $\left\vert h\right\vert $ sufficiently
small (e.g. for $|h|<|t|$), we have%
\begin{equation*}
\left\vert \left( t+h\right) ^{p}-t^{p}\right\vert =p\left\vert \sigma
\right\vert ^{p-1}\left\vert h\right\vert \leq p|2t|^{p-1}\left\vert
h\right\vert ,
\end{equation*}%
which implies 
\begin{equation*}
\Phi \left( T\left( t\right) ,T\left( t+h\right) \right) \leq \left\vert
h\right\vert \left\Vert \left\vert A\right\vert \right\Vert _{\mathcal{F}%
}\sum\limits_{p=1}^{\infty }\frac{1}{\left( p-1\right) !}\left(
|2t|\left\Vert \left\vert A\right\vert \right\Vert _{\mathcal{F}}\right)
^{p-1}.
\end{equation*}%
Passing to the limit as $h\rightarrow 0,$ it follows that $\underset{%
h\rightarrow 0^{+}}{\lim }\Phi \left( T\left( t+h\right) ,T\left( t\right)
\right) =0,$ i.e., $T\left( t\right) $ is continuous in $t\in \mathbf{R}.$

Now, let $t\geq 0$ and $h>0.$ Replacing (in the proof of relation $\left(
iii\right) $) $x$ by $T(t)(x)$, we get 
\begin{equation*}
\lim_{h\searrow 0}d(\frac{1}{h}\odot \lbrack T(h)(T(t)(x))\ominus
T(t)(x)],A(T(t)(x)))=0,
\end{equation*}%
and since $T(h)(T(t)(x))=T(t+h)(x)$, passing to supremum with $||x||_{%
\mathcal{F}}\leq 1$, it follows that%
\begin{equation*}
\underset{h\rightarrow 0^{+}}{\lim }\Phi \left( \frac{1}{h}\odot \left[
T\left( t+h\right) \ominus T\left( t\right) \right] ,A\left[ T\left(
t\right) \right] \right) =0.
\end{equation*}%
Similarly, denoting $t=t-h+h$, for $0<h<t,$ replacing in the proof of
relation (iii) $x$ by $T(t-h)(x)$ we obtain 
\begin{align*}
& d(\frac{1}{h}\odot \lbrack T(t)(x)\ominus T(t-h)(x)],A[T(t-h)(x)]) \\
& \leq ||T(t-h)(x)||_{\mathcal{F}}\frac{e^{h|||A|||_{\mathcal{F}%
}}-1-h|||A|||_{\mathcal{F}}}{h}.
\end{align*}%
Passing to the limit as $h\searrow 0$, from the continuity of $A$, $T(t)$
and reasoning as above, we arrive at 
\begin{equation*}
\lim_{h\searrow 0}d(\frac{1}{h}\odot \lbrack T(t)(x)\ominus
T(t-h)(x)],A(T(t)(x)))=0,
\end{equation*}%
and passing to supremum with $||x||_{\mathcal{F}}\leq 1$, it follows that 
\begin{equation*}
\underset{h\rightarrow 0^{+}}{\lim }\Phi \left( \frac{1}{h}\odot \left[
T\left( t\right) \ominus T\left( t-h\right) \right] ,A\left[ T\left(
t\right) \right] \right) =0.
\end{equation*}%
Now, let $t<0$ and $h>0$. Since by the above, (i) we have $%
S_{m}(k)(x)\rightarrow T(k)(x)$, for all $k<0$ and $x\in X$, repeating word
for word the proof in (iii), we immediately obtain 
\begin{equation*}
\lim_{k\nearrow 0}d(\frac{1}{k}\odot \lbrack T(k)(x)\ominus x],A(x))=0.
\end{equation*}%
Replacing $k$ by $-h$ and $x$ by $T(t)(x)$, by using (ii), we get 
\begin{equation*}
\lim_{h\searrow 0}d(-\frac{1}{h}\odot \lbrack T(t-h)(x)\ominus
T(t)(x)],A[T(t)(x)])=0.
\end{equation*}%
On the other hand, replacing $x$ by $T(t-k)(x)$ (where for $k$ sufficiently
close to zero, we have $t-k<0$) and then $k$ by $-h$, we arrive to 
\begin{equation*}
\lim_{h\searrow 0}d(-\frac{1}{h}\odot \lbrack T(t)(x)\ominus
T(t+h)(x)],A[T(t)(x)])=0,
\end{equation*}%
which proves (iv).

$\left( v\right) $ Let $t\geq 0$. We have%
\begin{equation*}
u^{^{\prime }}\left( t\right) =\left[ T\left( t\right) (u_{0})\right]
^{\prime }\oplus \left( \int\limits_{0}^{t}T\left( t-s\right) g\left(
s\right) ds\right) _{t}^{^{\prime }},
\end{equation*}%
where according to the point $\left( iv\right) ,$ we have $\left[ T\left(
t\right) (u_{0})\right] ^{^{\prime }}=A[T\left( t\right) (u_{0})].$ Denoting 
\begin{equation*}
F\left( t\right) =\int\limits_{0}^{t}T\left( t-s\right) (g\left( s\right)
)ds,
\end{equation*}%
we will show that $F^{^{\prime }}\left( t\right) $ exists and we shall find
it. To this end, let $h>0$. Denoting $F(t)=\int_{0}^{t}T(t-s)(g(s))ds$, we
have $t-s\geq 0$ and, using (ii), we get 
\begin{align*}
F(t+h)& =\int_{0}^{t+h}T(t-s+h)(g(s))ds=T(h)[\int_{0}^{t+h}T(t-s)(g(s))ds] \\
& =T(h)[F(t)\oplus \int_{t}^{t+h}T(t-s)(g(s))ds] \\
& =T(h)[F(t)]\oplus T(h)[\int_{t}^{t+h}T(t-s)(g(s))ds].
\end{align*}%
Using (iii), there exists $T(h)[F(t)]\ominus F(t)$, which implies 
\begin{equation*}
T(h)[F(t)]=[T(h)(F(t))\ominus F(t)]\oplus F(t).
\end{equation*}%
Replacing above and multiplying with $\frac{1}{h}\odot $, we obtain%
\begin{align*}
& \frac{1}{h}\odot \lbrack F(t+h)\ominus F(t)] \\
& =\frac{1}{h}\odot \lbrack T(h)(F(t))\ominus F(t)]\oplus T(h)[\frac{1}{h}%
\odot \int_{t}^{t+h}T(t-s)(g(s))ds].
\end{align*}%
Passing to the limit as $h\searrow 0,$ by (iii) it follows 
\begin{equation*}
\lim_{h\searrow 0}\frac{1}{h}\odot \lbrack F(t+h)\ominus F(t)]=A[F(t)]\oplus
\lim_{h\searrow 0}T(h)[\frac{1}{h}\odot \int_{t}^{t+h}T(t-s)(g(s))ds].
\end{equation*}%
Because for $h\searrow 0,$ we have $T(h)\rightarrow T(0)=I$, it remains to
show that 
\begin{equation*}
\lim_{h\searrow 0}\frac{1}{h}\odot \int_{t}^{t+h}T(t-s)(g(s))ds=g(t).
\end{equation*}%
Indeed, since $g(t)=T(0)(g(t))=\frac{1}{h}\odot \int_{t}^{t+h}T(0)(g(t))ds$,
we obtain 
\begin{equation*}
d(\frac{1}{h}\odot \int_{t}^{t+h}T(t-s)(g(s))ds,\frac{1}{h}\odot
\int_{t}^{t+h}T(0)(g(t))ds)\leq \frac{1}{h}\int_{t}^{t+h}H_{t}(s)ds,
\end{equation*}%
where $H_{t}(s)=d[T(t-s)(g(s)),T(0)(g(s))]$ is continuous on $[t,t+h]$ as
function of $s$, since $T\left( .\right) $ and $g$ are continuous.
Consequently,%
\begin{equation*}
\frac{1}{h}\int_{t}^{t+h}H_{t}(s)ds\rightarrow H_{t}(t)=0,\text{ as }%
h\searrow 0.
\end{equation*}%
Therefore, 
\begin{equation*}
\lim_{h\searrow 0}\frac{1}{h}\odot \lbrack F(t+h)\ominus F(t)]=A[F(t)]\oplus
g(t).
\end{equation*}%
On the other hand, for $0<h\leq t,$ we have $t-h\geq 0$ and $F(t)=F(u+h)$,
with $u=t-h\geq 0$. Reasoning as above, we obtain 
\begin{align*}
\frac{1}{h}\odot \lbrack F(u+h)\ominus F(u)]& =\frac{1}{h}\odot \lbrack
F(u+h)\ominus F(u)] \\
& =\frac{1}{h}\odot \lbrack T(h)[F(u)]\ominus F(u)]\oplus T(h)[\frac{1}{h}%
\odot \int_{u}^{u+h}T(u-s)(g(s))ds].
\end{align*}%
Because $F(t)$ is continuous (from the continuity of $T(.)$), clearly, $%
F(u)\rightarrow F(t)$ as $h\searrow 0$, which easily implies 
\begin{equation*}
\lim_{h\searrow 0}\frac{1}{h}\odot \lbrack T(h)(F(u))\ominus F(u)]=A[F(t)].
\end{equation*}%
Then, as in the case of the difference $\frac{1}{h}\odot \lbrack
F(t+h)\ominus F(t)]$, we obtain 
\begin{equation*}
\lim_{h\searrow 0}T(h)[\frac{1}{h}\odot \int_{u}^{u+h}T(u-s)(g(s))ds]=g(t).
\end{equation*}%
In conclusion, for $t\geq 0$ we have 
\begin{equation*}
\left( \int_{0}^{t}T(t-s)(g(s))ds\right) ^{\prime }=A[F(t)]\oplus g(t),
\end{equation*}%
which implies 
\begin{equation*}
u^{\prime }(t)=A[T(t)(u_{0})]\oplus A[F(t)]\oplus g(t)=A[u(t)]\oplus g(t).
\end{equation*}%
Now, let $t<0$ and $h>0$. As in the proof of the above point (iv), we use
the relation 
\begin{equation*}
\lim_{k\nearrow 0}d(\frac{1}{k}\odot \lbrack T(k)(x)\ominus x],A(x))=0,
\end{equation*}%
and repeating the above reasonings, we arrive at $u^{\prime
}(t)=A[u(t)]\oplus g(t)$. As before, the derivative $u^{\prime }(t)$ is
considered in the sense of Definition 2.7, (ii), 2). This proves (v) and
consequently the theorem. $\Box $

Theorem $3.9$ allows us to call $\left( e^{t\odot A}\right) _{t\geq 0}$ the
one-parameter fuzzy-semigroup generated by $A\in L\left( X\right) .$

\textbf{Remark.} Let $A,B:\mathbf{R}_{\mathcal{F}}\rightarrow \mathbf{R}_{%
\mathcal{F}}$ be defined by 
\begin{align*}
A(x)& =\left[ x_{-}\left( 1\right) -\int\limits_{0}^{1}x_{-}\left( r\right)
dr\right] \odot c, \\
B(x)& =\left[ x_{+}\left( 1\right) -\int\limits_{0}^{1}x_{+}\left( r\right)
dr\right] \odot c,
\end{align*}%
where $\left[ x\right] ^{r}=\left[ x_{-}\left( r\right) ,x_{+}\left(
r\right) \right] $ and $c\in \mathbf{R}_{\mathcal{F}}$ is a constant chosen
such that%
\begin{equation*}
\mu =c_{-}\left( 1\right) -\int\limits_{0}^{1}c_{-}\left( r\right) dr>0.
\end{equation*}%
By the Remark after Corollary $3.4$, we have $A,B\in L\left( \mathbf{R}_{%
\mathcal{F}}\right) .$ Then, by a simple calculation, we obtain that the
fuzzy-semigroups generated by $A$ and $B$ are 
\begin{equation*}
e^{t\odot A}(x)=x\oplus \frac{\left[ x_{-}\left( 1\right)
-\int\limits_{0}^{1}x_{-}\left( r\right) dr\right] }{\mu }\left( e^{t\mu
}-1\right) \odot c,t\geq 0
\end{equation*}%
and 
\begin{equation*}
e^{t\odot B}(x)x=x\oplus \frac{\left[ x_{+}\left( 0\right)
-\int\limits_{0}^{1}x_{+}\left( r\right) dr\right] }{\mu }\left( e^{t\mu
}-1\right) \odot c,t\geq 0,
\end{equation*}%
respectively. For example, by mathematical induction it easily follows 
\begin{equation*}
A^{n}(x)=\mu ^{n-1}\left( x_{-}\left( 1\right)
-\int\limits_{0}^{1}x_{-}\left( r\right) dr\right) \odot c
\end{equation*}%
and, using the Taylor series formula%
\begin{equation*}
e^{t\odot A}(x)=x\oplus \frac{t}{1!}\odot A(x)\oplus \frac{t^{2}}{2!}\odot
A^{2}(x)\oplus \frac{t^{3}}{3!}\odot A^{3}(x)\oplus ...\oplus \frac{t^{n}}{n!%
}\odot A^{n}(x)\oplus ...,
\end{equation*}%
we easily get the formula for $e^{t\odot A}(x).$

Another important result of this section is the following

\textbf{Theorem 3.10} \textit{For }$\mathit{(X,d)}$\textit{\ as in the
statement of Theorem 3.9 and $A\in L(X)$, let us define, formally,} 
\begin{equation*}
\mathit{T(t)=cosh[t\odot A]}
\end{equation*}%
in the sense that\textit{\ 
\begin{equation*}
\lim_{m\rightarrow +\infty }\Phi (T(t),\sum_{p=0}^{m}\frac{t^{2p}}{(2p)!}%
\odot A^{p})=0,
\end{equation*}%
(formally, we write $T(t)=cosh[t\odot A]=\sum_{p=0}^{+\infty }\frac{t^{2p}}{%
(2p)!}\odot A^{p}$). We have:}

\textit{(i) $T(t)\in L(X)$, for all $t\in $}$\mathbf{R}$\textit{; }

\textit{(ii) $T(t)$ is continuous as function of $t\in $}$\mathbf{R}$\textit{%
, $T(0)=I$ and if there exists $M>0$ such that $|||A|||_{\mathcal{F}%
}^{p}\leq M,\forall p=0,1,...,$ then T(t) is twice generalized
differentiable on }$\mathbf{R}$\textit{, $T^{\prime }(0)=\widetilde{0}_{X}$
and $T^{\prime \prime }(t)=A[T(t)]$. More exactly, for each $t\geq 0$ we
have 
\begin{equation*}
\lim_{h\searrow 0}\Phi (\frac{1}{h}\odot \lbrack T^{\prime }(t+h)\ominus
T^{\prime }(t)],A[T(t)])=0,
\end{equation*}%
\begin{equation*}
\lim_{h\searrow 0}\Phi (\frac{1}{h}\odot \lbrack T^{\prime }(t)\ominus
T^{\prime }(t-h)],A[T(t)])=0,
\end{equation*}%
where $T^{\prime }(t)$ is given by 
\begin{equation*}
\lim_{h\searrow 0}\Phi (\frac{1}{h}\odot \lbrack T(t+h)\ominus
T(t)],T^{\prime }(t))=0,
\end{equation*}%
\begin{equation*}
\lim_{h\searrow 0}\Phi (\frac{1}{h}\odot \lbrack T(t)\ominus
T(t-h)],T^{\prime }(t))=0,
\end{equation*}%
i.e., }%
\begin{equation*}
\mathit{T^{\prime }(t)(x)=sinh(t\odot A)=\sum_{p=1}^{+\infty }\frac{t^{2p-1}%
}{(2p-1)!}\odot A^{p}(x),}
\end{equation*}%
\textit{and for $t<0,$ we have 
\begin{equation*}
\lim_{h\searrow 0}\Phi (-\frac{1}{h}\odot \lbrack T^{\prime }(t)\ominus
T^{\prime }(t+h)],A[T(t)])=0,
\end{equation*}%
\begin{equation*}
\lim_{h\searrow 0}\Phi (-\frac{1}{h}\odot \lbrack T^{\prime }(t-h)\ominus
T^{\prime }(t)],A[T(t)])=0,
\end{equation*}%
where $T^{\prime }(t)$ is given by 
\begin{equation*}
\lim_{h\searrow 0}\Phi (-\frac{1}{h}\odot \lbrack T(t)\ominus
T(t+h)],T^{\prime }(t))=0,
\end{equation*}%
\begin{equation*}
\lim_{h\searrow 0}\Phi (-\frac{1}{h}\odot \lbrack T(t-h)\ominus
T(t)],T^{\prime }(t))=0.
\end{equation*}%
}

\textbf{Proof.} (i) Let us denote $C_{m}(t)(x)=\sum_{p=0}^{m}\frac{t^{2p}}{%
(2p)!}\odot A^{p}(x)$. Similarly, (see the proof of Theorem 3.9,$\left(
i\right) $) we obtain that $C_{m}(t)(x),m\in \mathbf{N}$ is a Cauchy
sequence and therefore its limit $T(t)=cosh(t\odot A)$ exists in $L(X)$.

(ii) The proof of continuity of $T(t)=cosh(t\odot A),$ $t\in \mathbf{R},$
i.e., $\lim_{h\rightarrow 0}\Phi (T(t),T(t+h))=0,\forall t\in \mathbf{R},$
is similar to the proof for $T(t)=e^{t\odot A},$ see the proof of Theorem
3.9,(iii). Now let $t\geq 0$ and $h>0$. We have 
\begin{align*}
C_{m}(t+h)(x)& =\sum_{p=0}^{m}\frac{(t+h)^{2p}}{(2p)!}\odot A^{p}(x) \\
& =C_{m}(t)(x)\oplus h\odot \lbrack t\odot A(x)\oplus \frac{t^{3}}{3!}\odot
A^{2}(x) \\
& \oplus ...\oplus \frac{t^{2m-1}}{(2m-1)!}\odot A^{m}(x)]\oplus
E_{m}(t,h,x),
\end{align*}%
where 
\begin{align*}
E_{m}(t,h,x)& =\frac{h^{2}}{2!}\odot A(x)\oplus \left[ \frac{{\binom{4}{2}}}{%
4!}t^{2}h^{2}+\frac{{\binom{4}{3}}}{4!}th^{3}+\frac{h^{4}}{4!}\right] \odot
A^{2}(x)\oplus ... \\
& \oplus \left[ \frac{{\binom{2k}{2}}}{(2k)!}t^{2k-2}h^{2}+\frac{{\binom{2k}{%
3}}}{(2k)!}t^{2k-3}h^{3}+...+\frac{h^{2k}}{(2k)!}\right] \odot
A^{k}(x)\oplus ... \\
& \oplus \left[ \frac{{\binom{2m}{2}}}{(2m)!}t^{2m-2}h^{2}+\frac{{\binom{2m}{%
3}}}{(2m)!}t^{2m-3}h^{3}+...+\frac{h^{2m}}{(2m)!}\right] \odot A^{m}(x).
\end{align*}

Let us denote 
\begin{equation*}
P_{m}(t)(x):=t\odot A(x)\oplus \frac{t^{3}}{3!}\odot A^{2}(x)\oplus
...\oplus \frac{t^{2m-1}}{(2m-1)!}\odot A^{m}(x).
\end{equation*}%
As in the proof of Theorem 3.9, we can show that $(P_{m}(t)(x))_{m}$ is a
Cauchy sequence, therefore it is convergent and let $\sinh (t\odot A)(x)$ be
its limit in $(X,\oplus ,\odot ,d)$. Moreover, note that for any $T>0$, $%
(C_{m}(t)(x))_{m}$ and $(P_{m}(t)(x))_{m}$ are uniformly Cauchy sequences on 
$[0,T]$. We get 
\begin{equation*}
\frac{1}{h}\odot \lbrack C_{m}(t+h)(x)\ominus C_{m}(t)(x)]=P_{m}(t)(x)\oplus 
\frac{1}{h}\odot E_{m}(t,h,x).
\end{equation*}%
On the other hand we have 
\begin{equation*}
d(\frac{1}{h}\odot E_{m}(t,h,x),\widetilde{0})\leq \frac{1}{h}%
M[ch_{m}(t+h)-ch_{m}(t)-hsh_{m}(t)],
\end{equation*}%
where $ch_{m}(t)=1+\frac{t^{2}}{2!}+...+\frac{t^{2m}}{(2m)!}$ and $%
sh_{m}(t)=t+\frac{t^{3}}{3!}+...+\frac{t^{2m-1}}{(2m-1)!}.$ Passing to the
limit as $h\searrow 0$, by the L'H\^{o}spital's rule, we obtain 
\begin{equation*}
\lim_{h\searrow 0}\frac{1}{h}[ch_{m}(t+h)-ch_{m}(t)-hsh_{m}(t)]=0,
\end{equation*}%
that is 
\begin{equation*}
\lim_{h\searrow 0}\frac{1}{h}\odot \lbrack C_{m}(t+h)(x)\ominus
C_{m}(t)(x)]=P_{m}(t)(x).
\end{equation*}%
Writing $C_{m}(t)(x)=C_{m}(t-h+h)(x)$, for $0<h<t$, by a similar reasoning,
we get 
\begin{equation*}
\lim_{h\searrow 0}\frac{1}{h}\odot \lbrack C_{m}(t)(x)\ominus
C_{m}(t-h)(x)]=P_{m}(t)(x).
\end{equation*}

Consequently, there exists the usual (Hukuhara) derivative in Definition
2.7, (i), $C_{m}^{\prime }(t)(x)=P_{m}(t)(x)$, for all $x\in X$. Thus, we
have obtained the following: for any $x\in X$, 
\begin{equation*}
\lim_{m\rightarrow +\infty }C_{m}(t)(x)=cosh(t\odot A)(x)
\end{equation*}%
and 
\begin{equation*}
\lim_{m\rightarrow +\infty }C_{m}^{\prime }(t)(x)=sinh(t\odot A)(x),
\end{equation*}%
uniformly with respect to $t\in \lbrack 0,T]$, for each $T>0$. Now, we show
that the sequences 
\begin{equation*}
F_{m}(h)=\frac{1}{h}\odot \lbrack C_{m}(t+h)(x)\ominus C_{m}(t)(x)],\text{ }%
m\in \mathbf{N}
\end{equation*}%
and 
\begin{equation*}
G_{m}(h)=\frac{1}{h}\odot \lbrack C_{m}(t)(x)\ominus C_{m}(t-h)(x)],\text{ }%
m\in \mathbf{N}
\end{equation*}%
are uniformly Cauchy sequences with respect to $h\in (0,T]$, for any $T>0$.
Let us prove the claim for $F_{m}(h)$ (the case of $G_{m}(h)$ is similar).
It is enough to show that $H_{m}(h)=\frac{1}{h}\odot E_{m}(t,h,x),m\in 
\mathbf{N},$ is a uniformly Cauchy sequence for $h\in (0,T],$ $T>0$
arbitrary (here $t>0$ and $x\in X$ are fixed). We obtain (by using the mean
value theorem) 
\begin{align*}
d(\frac{1}{h}\odot E_{m}(t,h,x),\widetilde{0})& \leq \frac{1}{h}%
M[hsh_{m}(\xi )\ominus hsh_{m}(t)] \\
& =M[sh_{m}(\xi )\ominus sh_{m}(t)] \\
& =M|\xi -t|ch_{m}(\eta )\leq Mhch_{m}(\eta ),
\end{align*}%
where $\xi \in (t,t+h),\eta \in (t,\xi ).$

When $h\in (0,T]$, we have $\xi \in (t,t+T]$ and $\eta \in (t,t+T]$ and $%
ch_{m}(\eta )\leq ch_{m}(t+T)\leq cosh(t+T)$, which implies 
\begin{equation*}
d(\frac{1}{h}\odot E_{m}(t,h,x),\widetilde{0})\leq Mhcosh(t+T),\forall m\in 
\mathbf{N},
\end{equation*}%
and 
\begin{equation*}
H_{m}(h)=\frac{1}{h}\odot E_{m}(t,h,x)\rightarrow \widetilde{0},
\end{equation*}%
uniformly with respect to $h\in (0,T]$ and, therefore, $(H_{m}(h))_{m}$ are
uniformly Cauchy for $h\in (0,T]$.

From the uniform convergence of the sequences $(C_{m}(t))_{m}$, $%
(C_{m}^{\prime }(t))_{m}$ with respect to $t\in (0,T]$ and of $%
(H_{m}(h))_{m} $ with respect to $h\in (0,T]$, by standard reasoning
(similar to that for the uniform convergence of derivative of the usual
sequences of functions), we obtain%
\begin{equation*}
T^{\prime }(t)(x)=[cosh(t\odot A)]^{\prime }(x)=sinh(t\odot A)(x).
\end{equation*}%
Repeating the procedure for $T^{\prime }(t)(x)$, we arrive at $T^{\prime
\prime }(t)=A[cosh(t\odot A)]=A[T(t)].$ The case $t<0$ is analagous. The
proof of the theorem is finished. $\Box $

\section{Applications to fuzzy differential equations}

\quad \quad In this section, we apply the main results of the previous
sections to solve fuzzy differential equations. As it was pointed out in the
Introduction, imprecision due to uncertainty or vagueness suggests of
considering fuzzy differential equations (i.e., whose solutions represents
fuzzy-number-valued functions), rather than random differential equations.
The simplest model is the following fuzzy Cauchy problem 
\begin{align}
\frac{du}{dt}\left( t\right) & =A[u\left( t\right) ],t\in I,  \tag{3} \\
u\left( t_{0}\right) & =u_{0},  \notag
\end{align}%
where $I$ is an interval, $u:I\rightarrow \mathbf{R}_{\mathcal{F}},u_{0}\in 
\mathbf{R}_{\mathcal{F}}$ and $A:C\left( I;\mathbf{R}_{\mathcal{F}}\right)
\rightarrow C\left( I;\mathbf{R}_{\mathcal{F}}\right) .$

The study of solutions of fuzzy differential equations has been considered
by e.g. $\left[ 3\right] -\left[ 7\right] ,\left[ 12\right] -\left[ 16\right]
,\left[ 18\right] .$ The main tools exploited are the so-called level set
and the differential inclusion methods. These methods are based on the idea
that to each fuzzy number $x\in \mathbf{R}_{\mathcal{F}}$ one can attach the
family of closed bounded real intervals%
\begin{equation*}
\left[ x\right] ^{r}=\left\{ t\in \mathbf{R;x}\left( t\right) \geq r\right\}
=\left[ x_{-}\left( r\right) ,x_{+}\left( r\right) \right] ,\text{ }r\in %
\left[ 0,1\right] .
\end{equation*}
However, none of these approaches use the powerful theory of semigroups of
operators to solve fuzzy differential equations, simply because the theory
was not developed until the present paper.

In what follows, we apply the theory of fuzzy-semigroups on some examples.
First, let us consider the general fuzzy Cauchy problem $\left( 3\right) .$
Here $\frac{du}{dt}$ means the derivative in the sense of the Definition $%
2.7,$ where $u\left( t\right) \in X,$ $t\geq 0,$ $A\in L\left( X\right) $
and $\left( X,d\right) $ can be chosen any of the following spaces: $\left( 
\mathbf{R}_{\mathcal{F}},D\right) ,$ $\left( l_{\mathbf{R}_{\mathcal{F}%
}}^{p},\rho _{p}\right) ,$ $\left( m_{\mathbf{R}_{\mathcal{F}}},\mu \right)
, $ $\left( c_{\mathbf{R}_{\mathcal{F}}},\mu \right) ,$ $\left( c_{\mathbf{R}%
_{\mathcal{F}}}^{\widetilde{0}},\mu \right) ,$ $\left( L^{p}\left( \left[ a,b%
\right] ;\mathbf{R}_{\mathcal{F}}\right) ,D_{p}\right) ,$ for $1\leq
p<\infty $, $\left( C\left( \left[ a,b\right] ;\mathbf{R}_{\mathcal{F}%
}\right) ,D^{\ast }\right) ,$ $\left( C^{p}\left( \left[ a,b\right] ;\mathbf{%
R}_{\mathcal{F}}\right) ,D_{p}^{\ast }\right) ,$ $p\in \mathbb{N},$ or any
finite Cartesian product of them endowed with the box metric. According to
Theorem $3.9,$ 
\begin{equation*}
u\left( t\right) =T\left( t\right) (u_{0})=e^{t\odot
A}(u_{0})=\sum\limits_{k=0}^{\infty }\frac{t^{k}}{k!}A^{k}(u_{0}),
\end{equation*}%
where $\sum $ is the sum with respect to $\oplus ,$ formally furnishes a
solution to the fuzzy Cauchy problem (3). We can first apply Theorem 3.9 to
the following fuzzy partial differential equation with initial conditions%
\begin{align}
& \left\{ 
\begin{array}{l}
\frac{du}{dt}\left( t,x\right) =v\left( t,x\right) \\ 
\frac{dv}{dt}\left( t,x\right) =u\left( t,x\right) ,t\geq 0,x\in \left[ a,b%
\right] ,%
\end{array}%
\right.  \tag{4} \\
& \left\{ 
\begin{array}{l}
u\left( 0,x\right) =u_{0}\left( x\right) , \\ 
v\left( 0,x\right) =v_{0}\left( x\right) ,x\in \left[ a,b\right] ,%
\end{array}%
\right.  \notag
\end{align}%
where $\frac{du}{dt}$ means the derivative of $u$ with respect to $t,$ see
Definition 2.7, (i), and $u\left( t,\cdot \right) ,$ $v(t,\cdot )\in C\left( %
\left[ a,b\right] ;\mathbf{R}_{\mathcal{F}}\right) ,$ for all $t\geq 0.$

Let us now set 
\begin{equation*}
A=\left( 
\begin{array}{ll}
0 & 1 \\ 
1 & 0%
\end{array}%
\right) ,w=\left( 
\begin{array}{l}
u \\ 
v%
\end{array}%
\right) \text{ and }w_{0}=\left( 
\begin{array}{l}
u_{0} \\ 
v_{0}%
\end{array}%
\right) .
\end{equation*}%
Then, problem (4) can be written as 
\begin{equation*}
\left\{ 
\begin{array}{l}
\frac{dw}{dt}\left( t\right) =\widetilde{A}[w\left( t\right) ],t\geq 0, \\ 
w\left( 0\right) =w_{0},%
\end{array}%
\right.
\end{equation*}%
where 
\begin{equation*}
\widetilde{A}\in \left[ C\left( \left[ a,b\right] ;\mathbf{R}_{\mathcal{F}%
}\right) \right] ^{2}\rightarrow \left[ C\left( \left[ a,b\right] ;\mathbf{R}%
_{\mathcal{F}}\right) \right] ^{2}
\end{equation*}%
is defined by%
\begin{equation*}
\widetilde{A}(w)=A\odot w=\left( 
\begin{array}{l}
v \\ 
u%
\end{array}%
\right) .
\end{equation*}%
Clearly, $A$ is a linear operator and continuous at each $w.$ Then, as in
e.g. $\left[ 9\right] ,$ we easily get that 
\begin{equation*}
T\left( t\right) =e^{t\odot \widetilde{A}}=\left( 
\begin{array}{ll}
\cosh \left( t\right) & \sinh \left( t\right) \\ 
\sinh \left( t\right) & \cosh \left( t\right)%
\end{array}%
\right) ,
\end{equation*}%
where $\sin h\left( t\right) =\frac{e^{t}-e^{-t}}{2},\cos h\left( t\right) =%
\frac{e^{t}+e^{-t}}{2},$ $t\geq 0.$ Consequently, a solution of $\left(
4\right) $ will be given by $w\left( t\right) =T\left( t\right) (w_{0}),$
i.e., 
\begin{equation*}
u\left( x,t\right) =\cosh \left( t\right) \odot u_{0}\left( x\right) \oplus
\sinh \left( t\right) \odot v_{0}\left( x\right) ,
\end{equation*}%
\begin{equation*}
v\left( x,t\right) =\sinh \left( t\right) \odot u_{0}\left( x\right) \oplus
\cosh \left( t\right) \odot v_{0}\left( x\right) ,
\end{equation*}%
such that $\left( u,v\right) $ is a solution of the fuzzy system $\left(
4\right) .$

Another example of a system of fuzzy differential equation is the following 
\begin{align}
& \left\{ 
\begin{array}{l}
\frac{du}{dt}\left( t,x\right) =u\left( t,x\right) \oplus v\left( t,x\right)
\\ 
\frac{dv}{dt}\left( t,x\right) =\left( -1\right) \odot u\left( t,x\right)
\oplus \left( -1\right) \odot v\left( t,x\right) ,t\geq 0,x\in \left[ a,b%
\right] ,%
\end{array}%
\right.  \tag{5} \\
& \left\{ 
\begin{array}{l}
u\left( 0,x\right) =u_{0}\left( x\right) , \\ 
v\left( 0,x\right) =v_{0}\left( x\right) ,x\in \left[ a,b\right] ,%
\end{array}%
\right.  \notag
\end{align}%
where $u\left( t,\cdot \right) ,$ $v\left( t,\cdot \right) \in C\left( \left[
a,b\right] ;\mathbf{R}_{\mathcal{F}}\right) ,$ for all $t\geq 0$ and the
derivatives are meant in the sense of Definition 2.7, (i). Let us set 
\begin{equation*}
A=\left( 
\begin{array}{ll}
1 & 1 \\ 
-1 & -1%
\end{array}%
\right) ,w=\left( 
\begin{array}{l}
u \\ 
v%
\end{array}%
\right) \text{ and }w_{0}=\left( 
\begin{array}{l}
u_{0} \\ 
v_{0}%
\end{array}%
\right) .
\end{equation*}
Then, problem (5)\ can be written as 
\begin{equation*}
\left\{ 
\begin{array}{l}
\frac{dw}{dt}\left( t\right) =\widetilde{A}[w\left( t\right) ],t\geq 0, \\ 
w\left( 0\right) =w_{0},%
\end{array}%
\right.
\end{equation*}%
where 
\begin{equation*}
\widetilde{A}\in \left[ C\left( \left[ a,b\right] ;\mathbf{R}_{\mathcal{F}%
}\right) \right] ^{2}\rightarrow \left[ C\left( \left[ a,b\right] ;\mathbf{R}%
_{\mathcal{F}}\right) \right] ^{2}
\end{equation*}%
is defined by 
\begin{equation*}
\widetilde{A}[w]=\left( 
\begin{array}{l}
u\oplus v \\ 
\left( -1\right) \odot \left[ u\oplus v\right]%
\end{array}%
\right) .
\end{equation*}%
Once more $A$ is a linear operator and continuous at each $w.$ We will
calculate $e^{t\odot \widetilde{A}}.$ In this case, we observe that 
\begin{equation*}
\widetilde{A}^{2}(w)=\left( 
\begin{array}{l}
u\oplus v\oplus \left( -1\right) \odot \left[ u\oplus v\right] \\ 
u\oplus v\oplus \left( -1\right) \odot \left[ u\oplus v\right]%
\end{array}%
\right) ,
\end{equation*}%
\begin{equation*}
\widetilde{A}^{n}(w)=2^{n-2}\odot \left( 
\begin{array}{l}
u\oplus v\oplus \left( -1\right) \odot \left[ u\oplus v\right] \\ 
u\oplus v\oplus \left( -1\right) \odot \left[ u\oplus v\right]%
\end{array}%
\right) ,n\geq 2,
\end{equation*}%
since denoting $E(u,v)=(u\oplus v)\oplus (-1)\odot (u\oplus v)$, we easily
get $(-1)\odot E(u,v)=E(u,v).$ Thus, 
\begin{align*}
T\left( t\right) (w)& =e^{t\odot \widetilde{A}}(w)=w\oplus \frac{t}{1!}\odot 
\widetilde{A}(w)\oplus \frac{t^{2}}{2!}\odot \widetilde{A}^{2}(w)\oplus
...\oplus \frac{t^{n}}{n!}\odot \widetilde{A}^{n}(w)\oplus ... \\
& =\left( 
\begin{array}{l}
u \\ 
v%
\end{array}%
\right) \oplus \left( 
\begin{array}{l}
t\odot \left( u\oplus v\right) \\ 
\left( -t\right) \odot \left( u\oplus v\right)%
\end{array}%
\right) \oplus \left( 
\begin{array}{l}
u\oplus v\oplus \left( -1\right) \odot \left[ u\oplus v\right] \\ 
u\oplus v\oplus \left( -1\right) \odot \left[ u\oplus v\right]%
\end{array}%
\right) \\
& \odot \left( \frac{t^{2}}{2!}+2\frac{t^{3}}{3!}+...+\frac{t^{n}}{n!}%
2^{n-2}+...\right) .
\end{align*}%
Let us denote 
\begin{equation*}
h\left( t\right) :=\sum\limits_{2}^{\infty }\frac{t^{k}}{k!}2^{k-2}.
\end{equation*}%
It easy to check that%
\begin{equation*}
\frac{1}{4}e^{2t}=\frac{1}{4}\sum\limits_{0}^{\infty }\frac{t^{k}}{k!}\left(
e^{2t}\right) _{t=0}^{\left( k\right) }=\frac{1}{4}\sum\limits_{0}^{\infty }%
\frac{t^{k}}{k!}2^{k}=h\left( t\right) +\frac{1}{4}\left( \frac{2t}{1!}%
+1\right) ,
\end{equation*}%
i.e., $h\left( t\right) =\frac{1}{4}\left( e^{2t}-2t-1\right) .$

Henceforth, we have%
\begin{equation*}
T\left( t\right) [w_{0}]=\left( 
\begin{array}{l}
u_{0} \\ 
v_{0}%
\end{array}%
\right) \oplus \left( 
\begin{array}{l}
t\odot \left( u_{0}\oplus v_{0}\right) \\ 
\left( -t\right) \odot \left( u_{0}\oplus v_{0}\right)%
\end{array}%
\right) \newline
\oplus \left( 
\begin{array}{l}
E\left( u_{0},v_{0}\right) \\ 
E\left( u_{0},v_{0}\right)%
\end{array}%
\right) \odot h\left( t\right) .
\end{equation*}%
Thus, an explicit solution of $(5)$ is given by 
\begin{equation*}
\left\{ 
\begin{array}{l}
u\left( t,x\right) =u_{0}\left( x\right) \oplus \left( t\right) \odot \left(
u_{0}\left( x\right) \oplus v_{0}\left( x\right) \right) \oplus E\left(
u_{0},v_{0}\right) \odot h\left( t\right) , \\ 
v\left( t,x\right) =v_{0}\left( x\right) \oplus \left( -t\right) \odot
\left( u_{0}\left( x\right) \oplus v_{0}\left( x\right) \right) \oplus
E\left( u_{0},v_{0}\right) \odot h\left( t\right) .%
\end{array}%
\right.
\end{equation*}

\textbf{Remark.} We don't know if the solution of the general fuzzy Cauchy
problem $(3)$ is unique. In the classical case, i.e., when the differential
equations and their corresponding solutions are in the setting of Banach
space valued functions, it is known that the uniqueness phenomenon holds.
For example, in this case, we know that the unique solution of (5) is given
by%
\begin{equation*}
u(t,x)=u_{0}(x)+t[u_{0}(x)+v_{0}(x)],\text{ }%
v(t,x)=v_{0}(x)-t[u_{0}(x)+v_{0}(x)]
\end{equation*}%
(note that it is obvious that $E(u_{0},v_{0})=0$). However, if $%
u_{0}(x),v_{0}(x)\in \mathbf{R}_{\mathcal{F}}\setminus \mathbf{R}$, then 
\begin{equation*}
u(t,x)=u_{0}(x)\oplus t\odot \lbrack u_{0}(x)\oplus v_{0}(x)],
\end{equation*}%
\begin{equation*}
v(t,x)=v_{0}(x)\oplus (-t)\odot \lbrack u_{0}(x)\oplus v_{0}(x)]
\end{equation*}%
is actually not a solution of (5), because%
\begin{equation*}
\frac{\partial u}{\partial t}=u_{0}(x)\oplus v_{0}(x)
\end{equation*}%
is a quantity which is essentially different from 
\begin{equation*}
u(t,x)\oplus v(t,x)=u_{0}(x)\oplus v_{0}(x)\oplus t\odot E(u_{0},v_{0}).
\end{equation*}

In what follows we apply Theorem 3.10 on two more examples. Let us consider
the following system of fuzzy partial differential equation with initial
conditions 
\begin{align}
& \left\{ 
\begin{array}{l}
\frac{{\partial u}^{2}}{\partial {t}^{2}}\left( t,x\right) =u\left(
t,x\right) \oplus v\left( t,x\right) \\ 
\frac{{\partial v}^{2}}{\partial {t}^{2}}\left( t,x\right) =\left( -1\right)
\odot u\left( t,x\right) \oplus \left( -1\right) \odot v\left( t,x\right)
,t\geq 0,x\in \left[ a,b\right] ,%
\end{array}%
\right.  \tag{6} \\
& \left\{ 
\begin{array}{l}
u\left( 0,x\right) =u_{0}\left( x\right) , \\ 
v\left( 0,x\right) =v_{0}\left( x\right) ,\frac{\partial u}{\partial t}(0,x)=%
\frac{\partial v}{\partial t}(0,x)=\widetilde{0},x\in \left[ a,b\right] ,%
\end{array}%
\right.  \notag
\end{align}%
where $u\left( t,\cdot \right) ,$ $v\left( t,\cdot \right) \in C\left( \left[
a,b\right] ;\mathbf{R}_{\mathcal{F}}\right) ,$ for all $t\geq 0$ and where
the second derivatives $\frac{{\partial u}^{2}}{\partial {t}^{2}},$ $\frac{{%
\partial v}^{2}}{\partial {t}^{2}}$ are meant in the sense of Definition
2.7, (i). As in the previous example, let us set%
\begin{equation*}
A=\left( 
\begin{array}{ll}
1 & 1 \\ 
-1 & -1%
\end{array}%
\right) ,w=\left( 
\begin{array}{l}
u \\ 
v%
\end{array}%
\right) ,\text{ }\widetilde{0}=\left( \widetilde{0},\widetilde{0}\right) 
\text{ and }w_{0}=\left( 
\begin{array}{l}
u_{0} \\ 
v_{0}%
\end{array}%
\right) .
\end{equation*}%
Then, system (6) can be written as 
\begin{equation*}
\left\{ 
\begin{array}{l}
\frac{{\partial }^{2}w(t)}{\partial {t}^{2}}\left( t\right) =\widetilde{A}%
[w\left( t\right) ],t\geq 0, \\ 
w\left( 0\right) =w_{0},\frac{\partial w}{\partial t}(0)=\widetilde{0},%
\end{array}%
\right.
\end{equation*}%
where 
\begin{equation*}
\widetilde{A}\in \left[ C\left( \left[ a,b\right] ;\mathbf{R}_{\mathcal{F}%
}\right) \right] ^{2}\rightarrow \left[ C\left( \left[ a,b\right] ;\mathbf{R}%
_{\mathcal{F}}\right) \right] ^{2}
\end{equation*}%
is defined by 
\begin{equation*}
\widetilde{A}[w]=\left( 
\begin{array}{l}
u\oplus v \\ 
\left( -1\right) \odot \left[ u\oplus v\right]%
\end{array}%
\right) .
\end{equation*}%
According to Theorem 3.10, a solution is given by the formula 
\begin{equation*}
T(t)(w_{0})=\cosh (t\odot \widetilde{A})(w_{0}).
\end{equation*}%
We will calculate $\cosh (t\odot \widetilde{A})$. Similarly to the previous
example, we have 
\begin{equation*}
\widetilde{A}^{2}(w)=\left( 
\begin{array}{l}
u\oplus v\oplus \left( -1\right) \odot \left[ u\oplus v\right] \\ 
u\oplus v\oplus \left( -1\right) \odot \left[ u\oplus v\right]%
\end{array}%
\right) ,
\end{equation*}%
\begin{equation*}
\widetilde{A}^{n}(w)=2^{n-2}\odot \left( 
\begin{array}{l}
u\oplus v\oplus \left( -1\right) \odot \left[ u\oplus v\right] \\ 
u\oplus v\oplus \left( -1\right) \odot \left[ u\oplus v\right]%
\end{array}%
\right) ,\text{ }n\geq 2,
\end{equation*}%
so that%
\begin{align*}
T\left( t\right) (w)& =\cosh (t\odot \widetilde{A})(w)=w\oplus \frac{t^{2}}{%
2!}\odot \widetilde{A}(w)\oplus \frac{t^{4}}{4!}\odot \widetilde{A}^{2}(w) \\
& \oplus ...\oplus \frac{t^{2n}}{(2n)!}\odot \widetilde{A}^{n}(w)\oplus ...
\\
& =\left( 
\begin{array}{l}
u \\ 
v%
\end{array}%
\right) \oplus \left( 
\begin{array}{l}
\frac{t^{2}}{2!}\odot \left( u\oplus v\right) \\ 
\left( -\frac{t^{2}}{2!}\right) \odot \left( u\oplus v\right)%
\end{array}%
\right) \\
& \oplus \left( 
\begin{array}{l}
u\oplus v\oplus \left( -1\right) \odot \left[ u\oplus v\right] \\ 
u\oplus v\oplus \left( -1\right) \odot \left[ u\oplus v\right]%
\end{array}%
\right) \odot \left( \frac{t^{4}}{4!}+2\frac{t^{6}}{6!}+...+\frac{t^{2n}}{%
(2n)!}2^{n-2}+...\right) .
\end{align*}%
Let us set 
\begin{equation*}
h\left( t\right) :=\sum\limits_{2}^{\infty }\frac{t^{2k}}{(2k)!}2^{k-2}
\end{equation*}%
and 
\begin{equation*}
E\left( u,v\right) :=u\oplus v\oplus \left( -1\right) \odot \left[ u\oplus v%
\right] .
\end{equation*}%
It easy to check that 
\begin{equation*}
h(t)=\frac{1}{4}\sum\limits_{2}^{\infty }\frac{(t\sqrt{2})^{2k}}{(2k)!}=%
\frac{1}{4}[\cosh (t\sqrt{2})-t^{2}-1].
\end{equation*}%
Finally, we then have 
\begin{equation*}
T\left( t\right) [w_{0}]=\left( 
\begin{array}{l}
u_{0} \\ 
v_{0}%
\end{array}%
\right) \oplus \left( 
\begin{array}{l}
\frac{t^{2}}{2}\odot \left( u_{0}\oplus v_{0}\right) \\ 
\left( -\frac{t^{2}}{2}\right) \odot \left( u_{0}\oplus v_{0}\right)%
\end{array}%
\right) \newline
\oplus \left( 
\begin{array}{l}
E\left( u_{0},v_{0}\right) \\ 
E\left( u_{0},v_{0}\right)%
\end{array}%
\right) \odot h\left( t\right) ,
\end{equation*}%
such that an explicit solution of $(6)$ is given by 
\begin{equation*}
\left\{ 
\begin{array}{l}
u\left( t,x\right) =u_{0}\left( x\right) \oplus \left( \frac{t^{2}}{2}%
\right) \odot \left( u_{0}\left( x\right) \oplus v_{0}\left( x\right)
\right) \oplus E\left( u_{0},v_{0}\right) \odot h\left( t\right) , \\ 
v\left( t,x\right) =v_{0}\left( x\right) \oplus \left( -\frac{t^{2}}{2}%
\right) \odot \left( u_{0}\left( x\right) \oplus v_{0}\left( x\right)
\right) \oplus E\left( u_{0},v_{0}\right) \odot h\left( t\right) .%
\end{array}%
\right.
\end{equation*}

\textbf{Remark.} As in the previous example we don't know if this solution
is unique, but in the case when $u_{0}(x)$ and $v_{0}(x)$ are real-valued
functions, we obtain the unique solution (since $E(u_{0},v_{0})=0$) 
\begin{equation*}
u(t,x)=u_{0}(x)+\frac{t^{2}}{2}[u_{0}(x)+v_{0}(x)],\text{ }v(t,x)=v_{0}(x)-%
\frac{t^{2}}{2}[u_{0}(x)+v_{0}(x)].
\end{equation*}%
Note that if $u_{0}(x),v_{0}(x)\in \mathbf{R}_{\mathcal{F}}\backslash 
\mathbf{R}$, then 
\begin{equation*}
u(t,x)=u_{0}(x)\oplus \frac{t^{2}}{2}\odot \lbrack u_{0}(x)\oplus v_{0}(x)],%
\text{ }v(t,x)=v_{0}(x)\oplus (-\frac{t^{2}}{2})\odot \lbrack u_{0}(x)\oplus
v_{0}(x)]
\end{equation*}%
is actually not a solution of (6).

As our final example, let us consider the initial value problem for the
fuzzy wave equation 
\begin{equation*}
\left\{ 
\begin{array}{l}
\frac{{\partial u}^{2}}{\partial {t}^{2}}\left( t,x\right) =\frac{{\partial u%
}^{2}}{\partial {x}^{2}}\left( t,x\right) \\ 
u\left( 0,x\right) =u_{1}\left( x\right) , \\ 
\frac{\partial u}{\partial t}(0,x)=u_{2}(x),x\in \left[ a,b\right] ,%
\end{array}%
\right.
\end{equation*}%
where $\frac{\partial u}{\partial t},\frac{{\partial }^{2}u}{\partial t^{2}}$
are Hukuhara-kind derivatives and $u_{1}:\mathbf{R}\rightarrow X$ is
supposed to be Hukuhara differentiable of any order (in the sense of
Definition 2.7, (i)), such that there exists $M>0$ satisfying 
\begin{equation*}
d(\widetilde{0},u_{1}^{(2p)}(x))\leq M,\forall x\in X,p=0,1,...,.
\end{equation*}%
Then $A(u)=\frac{{\partial }^{2}u}{\partial x^{2}}$ and%
\begin{equation*}
\cosh (t\odot A)(u_{1}(x))=\sum_{p=0}^{+\infty }\frac{t^{2p}}{(2p)!}\odot
u_{1}^{(2p)}(x)
\end{equation*}%
is a convergent series in $(X,d)$, with $\sum_{p=0}^{m}\frac{t^{2p}}{(2p)!}%
\odot u_{1}^{(2p)}(x)$ being a uniform approximation (with respect to $t$ on
compact subintervals) of $\cosh (t\odot A)(u_{1}(x))$. As a consequence, 
\begin{equation*}
u(t,x)=\cosh (t\odot A)\oplus t\odot u_{2}(x)
\end{equation*}%
furnishes a solution to the above fuzzy wave equation. Furthermore,%
\begin{equation*}
u_{m}(t,x)=\sum_{p=0}^{m}\frac{t^{2p}}{(2p)!}\odot u_{1}^{(2p)}(x)\oplus
t\odot u_{2}(x),\text{ }m\in \mathbf{N},
\end{equation*}%
represents a sequence of uniform approximation (with respect to $t$ on
compact subintervals) for the exact solution $u(t,x)$.

\textbf{Remark.} In view of Theorems 3.9 and 3.10 we can also take into the
above examples the differentiability in the generalized sense as well (see
Definition 2.7, (ii)).

We can conclude that, as in the classical theory of differential equations,
the theory of semigroups of operators developed in this paper may become a
powerful tool to study solutions of fuzzy partial differential equations.

\end{document}